\documentclass[graybox]{svmult}

\pdfoutput=1

\usepackage{amsmath}									
\usepackage{amsbsy}
\usepackage{color}
\usepackage{esint}
\usepackage{comment}
\usepackage{soul}

\usepackage[T1]{fontenc}
\usepackage[bitstream-charter]{mathdesign}
\usepackage[latin1]{inputenc}	

\usepackage{mathptmx}       
\usepackage{helvet}         
\usepackage{courier}        
\usepackage{type1cm}        
\usepackage{makeidx}         
\usepackage{graphicx}        
\usepackage{multicol}        
\usepackage[bottom]{footmisc}

\newcommand*{\Var}{\operatorname{Var}}

\newcommand*{\wt}{\widetilde}

\newcommand*{\EE}{\mathbb E}
\newcommand*{\PP}{\mathbb P}

\newcommand*{\RR}{\mathbb R}
\newcommand*{\NN}{\mathbb N}

\newcommand*{\bbN}{\mathbb N}

\newcommand*{\bbR}{\mathbb R}

\newcommand*{\cF}{\mathcal F}

\newcommand*{\cJ}{\mathcal J}

\DeclareMathOperator{\sgn}{sgn}
\DeclareMathOperator{\arsinh}{arsinh}
\DeclareMathOperator{\sech}{sech}

\renewcommand*{\doteq}{:=}

\makeindex             

\begin{document}

\title*{Regression-based variance reduction approach for strong approximation schemes}
\author{Denis Belomestny, Stefan H\"afner and Mikhail Urusov}
\institute{Denis Belomestny \at Duisburg-Essen University, Essen, Germany, \email{denis.belomestny@uni-due.de}\and Stefan H\"afner \at PricewaterhouseCoopers GmbH, Frankfurt, Germany, \email{stefan.haefner@de.pwc.com}\and Mikhail Urusov \at Duisburg-Essen University, Essen, Germany, \email{mikhail.urusov@uni-due.de}}
%
%
\maketitle

\abstract{In this paper we present a novel approach towards  variance reduction for discretised diffusion processes. The proposed approach involves specially constructed control variates and allows for a significant reduction in the variance for the terminal functionals. In this way the complexity order of the standard Monte Carlo algorithm ($\varepsilon^{-3}$) can be reduced down to $\varepsilon^{-2}\sqrt{\left|\log(\varepsilon)\right|}$ in case of the Euler scheme with $\varepsilon$ being the precision to be achieved. These theoretical results are illustrated by several numerical examples.}

\section{Introduction}
\label{sec:1}

Let \(T>0\) be a fixed time horizon.
Consider a $d$-dimensional diffusion process
$(X_t)_{t\in[0,T]}$
defined on a filtered probability space
$(\Omega,\cF,(\cF_t)_{t\in[0,T]},\PP)$
by the It\^o stochastic differential equation
\begin{align}
\label{x_sde}
dX_t
=\mu(X_t)\,dt
+\sigma(X_t)\,dW_{t},
\quad X_{0}=x_0\in{\RR^d},
\end{align}
for Lipschitz continuous functions
\(\mu\colon\RR^d\to\RR^d\)
and
\(\sigma\colon\RR^d
\to\RR^{d\times m}\),
where \((W_t)_{t\in[0,T]}\)
is a standard \(m\)-dimensional
$(\cF_t)$-Brownian motion.
Suppose we want to find a continuous function 
\begin{align*}
u=u(t,x)\colon[0,T]\times\RR^d\to\RR,
\end{align*}
which has a continuous first derivative with respect to \(t\) and continuous first and second derivatives with respect to the components of \(x\) on $[0,T)\times\RR^d$,
such that it solves the partial differential equation 
\begin{align}
\label{eq:Cauchy_prob}
\frac{\partial u}{\partial t}+\mathcal{L} u&=0 \quad \mbox{ on }  [0,T)\times \RR^d,\\
\label{eq:term_cond}
u(T,x)&=f(x)\quad 	\mbox{ for } x\in \RR^d,
\end{align}
where $f$ is a given Borel function on~$\RR^d$.
Here, \(\mathcal{L}\) is the differential operator associated with the equation \eqref{x_sde}:
\begin{align*}
(\mathcal{L}u)(t,x)
\doteq\sum_{k=1}^{d}\mu_{k}(x)\frac{\partial u}{\partial x_{k}}(t,x)
+\frac{1}{2}\sum_{k,l=1}^{d}(\sigma\sigma^\top)_{kl}(x)\frac{\partial^{2} u}{\partial x_{k}\partial x_{l}}(t,x),
\end{align*}
where $\sigma^\top$ denotes the transpose of~$\sigma$.
Under appropriate conditions on
$\mu$, $\sigma$ and~$f,$ there is a solution
of the Cauchy
problem \eqref{eq:Cauchy_prob}--\eqref{eq:term_cond},
which is unique in the class of solutions
satisfying certain growth conditions,
and it has the following Feynman-Kac
stochastic representation
\begin{align}
\label{eq:15022017a1}
u(t,x)=\EE[ f(X_{T}^{t,x})]
\end{align}
(see Section~5.7 in~\cite{karatzas2012brownian}),
where $X^{t,x}$ denotes the solution started
at time $t$ in point~$x$. Moreover it holds
\begin{align*}
\EE[f(X_{T}^{0,x})|X_{t}^{0,x}]=u(t,X_{t}^{0,x}), \quad \mbox{ a.s. }
\end{align*}
for \(t\in [0,T]\) and 
\begin{align}
\label{repr_contr_var}
f(X_{T}^{0,x})=\EE [f(X_{T}^{0,x})]+M^{*}_T, \quad \mbox{ a.s. }
\end{align}
with 
\begin{align}
\label{MT_cont}
M^*_T\doteq \int_{0}^{T}\nabla_x u (t,X_{t}^{0,x})\,\sigma(X_{t}^{0,x})\,dW_{t}
\equiv\int_{0}^{T}\sum_{k=1}^d \frac{\partial u}{\partial x_k} (t,X_{t}^{0,x})\sum_{i=1}^m\sigma_{ki}(X_{t}^{0,x})\,dW^i_{t}.
\end{align}

The standard Monte Carlo (SMC)
approach for computing \(u(0,x)\) at a fixed point \(x\in \bbR^d\) basically consists of three steps. First, an approximation \(\overline{X}_T\) for \(X^{0,x}_T\) is constructed via a time discretisation in equation~\eqref{x_sde} (we refer to~\cite{KP} for a nice overview of various discretisation schemes).
In this paper we focus on the Euler-Maruyama approximation
to the exact solution (the Euler scheme).
Next, \(N_0\) independent copies of the approximation
\(\overline{X}_T\)  are generated, and, finally, a Monte Carlo estimate \(V_{N_0}\) is defined as the average of the values of \(f\) at simulated points:
\begin{align}\label{eq:2512a1}
V_{N_0}\doteq\frac{1}{N_0}\sum_{n=1}^{N_0} f\Bigl(\overline{X}_T^{(n)}\Bigr).
\end{align} 
In the computation of $u(0,x)=\EE [f(X^{0,x}_T)]$
by the SMC approach
there are two types of error inherent:
the discretisation error
$\EE [f(X^{0,x}_T)]-\EE [f(\overline X_{T})]$
and the Monte Carlo (statistical) error,
which results from the substitution of
$\EE [f(\overline{X}_{T})]$
with the sample average~$V_{N_0}$.
The aim of variance reduction methods is to reduce the statistical error.  For example, in the so-called control variate variance reduction approach
one looks for a random variable
\(\xi\) with \(\EE \xi=0\),
which can be simulated,
such that
the variance of the difference
\(f(\overline{X}_{T})-\xi\) is minimised, that~is,
\begin{align*}
\Var[f(\overline{X}_{T})-\xi]\to\min\text{ under } \EE\xi=0.
\end{align*}

The use of control variates for solving \eqref{x_sde} via Monte Carlo path simulation approach was initiated by Newton~\cite{newton1994variance} and further developed in Milstein and Tretyakov~\cite{milstein2009practical}. In fact, the construction of the appropriate control variates  in the above two papers  essentially relies on identities
\eqref{repr_contr_var} and~\eqref{MT_cont}
implying that the zero-mean random variable $M^*_{T}$
can be viewed as an optimal control variate,
since 
\begin{align*}
\Var[f(X^{0,x}_{T})-M^*_{T}]
=\Var[\EE f(X^{0,x}_{T})]=0.
\end{align*}
Let us note that it would be desirable to have
a control variate  reducing the variance
of $f(\overline{X}_{T})$ rather than the one of $f(X^{0,x}_T)$
because we simulate from the distribution of $f(\overline{X}_{T})$
and not from the one of $f(X^{0,x}_T)$. Moreover, the control variate \(M^*_{T}\) cannot be directly computed, since the function \(u(t,x)\) is unknown. This is why Milstein and Tretyakov~\cite{milstein2009practical} proposed to use regression for getting a preliminary  approximation for  \(u(t,x)\) in a first step.

The contribution of our work is as follows.
We propose an approach for the construction of control variates that reduce the variance of \(f(\overline{X}_{T})\),
i.e.\ we perform variance reduction
not for the exact but rather for the discretised process.
A nice by-product is that our control variates
can be computed in a rather simple way,
and less assumptions are required in our case,
than one would require to construct control variates
based on on the exact solution.
Moreover, we present bounds for the regression error
involved in the construction of our control variates
and perform the complexity analysis
(these are not present in~\cite{milstein2009practical}),
which is also helpful for designing numerical experiments.
We are able to achieve a sufficient convergence order of the resulting variance, which in turn leads to a significant complexity reduction as compared to the SMC algorithm. Other examples of algorithms with this property include the analogous regression-based variance reduction approach for weak approximation schemes of~\cite{belomestny2016variance},
the multilevel Monte Carlo (MLMC) algorithm of~\cite{giles2008multilevel} and the quadrature-based algorithm of~\cite{muller2015complexity}. 
\par
Summing up, we propose a new regression-type approach for the construction of 
control variates in case of the Euler scheme.
It takes advantage of the smoothness in $\mu$, $\sigma$ and $f$ (which is needed for nice convergence properties of regression methods) in order to significantly reduce the variance of the random variable $f(\overline{X}_{T})$.

This work is organised as follows. In Section~\ref{sec:2}
we describe the construction of control variates for strong approximation schemes. Section~\ref{sec:3} describes the use of regression algorithms for the construction of control variates and analyses their convergence. A complexity analysis of the variance reduced Monte Carlo algorithm is conducted in Section~\ref{sec:4}. Section~\ref{sec:5} is devoted to a simulation study. Finally, all proofs are collected in Section~\ref{sec:proofs}.

\bigskip
\textbf{Notational convention.}
Throughout, elements of $\bbR^d$ (resp.~$\bbR^{1\times d}$)
are understood as column-vectors (resp.~row-vectors).
Generally, most vectors in what follows are column-vectors.
However, gradients of functions and some vectors defined
via them are row-vectors.
Finally, we record our standing assumption
that we do not repeat explicitly in the sequel.

\medskip\noindent
\textbf{Standing assumption.}
The coefficients $\mu$ and $\sigma$
in~\eqref{x_sde}
are globally Lipschitz functions.

\section{Control variates for strong approximation schemes}
\label{sec:2}
To begin with, we introduce some notations,
which will be frequently used in the sequel.
Throughout this paper,
$\bbN_0\doteq\bbN\cup\{0\}$
denotes the set of nonnegative integers,
$J\in\bbN$ denotes the time discretisation parameter,
we set $\Delta\doteq T/J$
and consider discretisation schemes
defined on the grid
$\{t_j=j\Delta:j=0,\ldots,J\}$. We set
$\Delta_j W\doteq W_{j\Delta}-W_{(j-1)\Delta}$, and by $W^i$ we denote the $i$-th component of the vector $W$. 
Further, for $k\in\bbN_0$,
$H_k\colon\bbR\to\bbR$
stands for the (normalised) $k$-th Hermite polynomial,~i.e.
\begin{align*}
  H_k(x) 
  \doteq 
  \frac{(-1)^k}
    {\sqrt{k!}}
  e^{\frac{x^2}{2}}
  \frac{d^k}{dx^k}e^{-\frac{x^2}{2}},
  \quad x\in\bbR.
\end{align*}
Notice that $H_0\equiv1$,
$H_1(x)=x$, $H_2(x)=\frac1{\sqrt{2}}(x^2-1)$.

\subsection{Series representation}
\label{series}
Let us consider a scheme,
where $d$-dimensional approximations
$X_{\Delta,j\Delta}$,
$j=0,\ldots,J$, satisfy $X_{\Delta,0}=x_0$ and
\begin{align}
\label{eq:defX}
    X_{\Delta, j\Delta} 
    = 
    \Phi_{\Delta}\left(
       X_{\Delta, (j-1)\Delta},  \Delta_{j} W
    \right),
\end{align}
where $\Delta_j W:=W_{j\Delta}-W_{(j-1)\Delta}$,
for some Borel measurable functions
$\Phi_{\Delta}\colon\bbR^{d\times m}\to\bbR^d$
(clearly, the Euler scheme is a special case of this setting).

\begin{theorem}\label{thm:ChaosDecompNum}
Let $f\colon \RR^d \rightarrow \RR$ be
a Borel measurable function
such that it holds $\EE | f(X_{\Delta, T}) |^2 < \infty$.
Then we have the representation
(cf.~Theorem~2.1 in~\cite{belomestny2016variance})
\begin{align}
\label{eq:chaos}
f(X_{\Delta,T})=\EE[f(X_{\Delta,T})]+  
\sum_{j=1}^{J}
\sum_{k\in\bbN_0^m\setminus\{0_m\}}
a_{j,k}(X_{\Delta, (j-1)\Delta})
\prod_{r=1}^m
H_{k_r}\left(\frac{\Delta_{j} W^{r}}{\sqrt{\Delta}}\right),
\end{align}
where $k=(k_1,\ldots,k_m)$ and $0_m\doteq(0,\ldots,0)\in\bbR^m$
(in the second summation),
and the coefficients
$a_{j,k}\colon\bbR^d\to\bbR$
are given by the formula
\begin{align}
\label{eq:wiener_pcv}
a_{j,k}(x)=\EE\left[ f(X_{\Delta,T})
\prod_{r=1}^m H_{k_r}
\left(\frac{\Delta_{j} W^{r}}{\sqrt{\Delta}}\right)
\bigg|\,X_{\Delta, (j-1)\Delta}=x \right],
\end{align}
for all $j\in \{1,\ldots,J\}$
and $k\in\bbN_0^m\setminus\{0_m\}$. 
\end{theorem}

\begin{remark}
Representation~\eqref{eq:chaos} shows that we have a perfect control variate, namely
\begin{eqnarray}
\label{cv_perfect}
M_{\Delta,T}\doteq
\sum_{j=1}^{J}
\sum_{k\in\bbN_0^m\setminus\{0_m\}}
a_{j,k}(X_{\Delta, (j-1)\Delta})
\prod_{r=1}^m
H_{k_r}\left(\frac{\Delta_{j} W^{r}}{\sqrt{\Delta}}\right),
\end{eqnarray}
for the functional \(f(X_{\Delta, T})\),
i.e.\ \(\Var[f(X_{\Delta, T})-M_{\Delta,T}]=0\).
\end{remark}

The control variate $M_{\Delta,T}$ is not implementable
because of the infinite summation in~\eqref{cv_perfect}
and because the coefficients $a_{j,k}$ are unknown.
In the later sections we estimate
the unknown coefficients in this and other
(related) representations via regression
and present bounds for the estimation error.

Now we introduce the following ``truncated'' control variate
\begin{align}\label{cv_truncated}
M^{{ser},1}_{\Delta,T}\doteq\sum_{j=1}^J\sum_{i=1}^ma_{j,e_i}(X_{\Delta,(j-1)\Delta})\frac{\Delta_j W^i}{\sqrt{\Delta}},
\end{align}
where $e_i$ denotes the $i$-th unit vector in $\bbR^m$.
The superscript ``ser'' comes from ``series''.
In the next subsection,
performing a quite different argumentation,
we derive another control variate,
which will turn out to be theoretically equivalent
to~$M^{{ser},1}_{\Delta,T}$.

\subsection{Integral representation}
\label{subs:int_repr}
\textbf{Integral representation for the exact solution.}
We first motivate what we call
``integral representation for the discretisation'',
which will be presented below in this subsection,
in that we recall in more detail
the main idea of constructing control variates
in Milstein and Tretyakov~\cite{milstein2009practical}.
As was already mentioned in the introduction,
the control variate in~\cite{milstein2009practical}
is an approximation of $M_T^*$ of~\eqref{MT_cont},
where the function $u$ is given in~\eqref{eq:15022017a1}
and is therefore unknown,
which rises the question about
a possible practical implementation
of~\eqref{MT_cont}.

To this end, let us define the ``derivative processes'' $\delta^iX^k_{s,x}(t)\doteq \frac{\partial X^k_{s,x}(t)}{\partial x_i}$ for $i,k\in\left\{1,\ldots,d\right\}$,
where $X^k_{s,x}(t)$ means the $k$-th component
of the solution of~\eqref{x_sde}
started at time~$s$ in~$x$
evaluated at time~$t\ge s$,
and simply write $\delta^i X^k_t$ rather than $\delta^i X^k_{0,x_0}(t)$ below. Further, we define the matrix $\delta X_t\doteq\begin{pmatrix}\delta^1X^1_t & \cdots & \delta^d X^1_t\\
\vdots & \ddots & \vdots\\
\delta^1X^d_t& \cdots & \delta^dX^d_t\end{pmatrix}\in\bbR^{d\times d}$ as well as the vectors $\delta^i X_t\doteq\begin{pmatrix}\delta^iX^1_t &  \cdots & \delta^iX^d_t\end{pmatrix}^\top\in\bbR^d$.
Assuming $\mu,\sigma\in C^1$,
we notice that $\delta^i X_t$ satisfies the following SDE
\begin{align}
\label{dx_sde}
d\delta^i X_t
=\sum_{k=1}^d\delta^i X^k_t\left[\frac{\partial\mu(X_t)}{\partial x_k }\,dt
+\frac{\partial\sigma(X_t)}{\partial x_k}\,dW_{t}\right],
\quad \delta^i X^k_{0}=\left\{\begin{array}{ll}1, & i=k\\ 0, & i\neq k\end{array}\right..
\end{align}
Milstein and Tretyakov~\cite{milstein2009practical}
exploit~\eqref{dx_sde} to prove that, provided $f,\mu,\sigma\in C^1$,
the integral in~\eqref{MT_cont} can be expressed by means of $\delta X_t$ as follows
\begin{align}
\label{cv_milstein}
M^*_T\doteq \int\limits_{0}^{T}\nabla_xu(t,X_{t})\,\sigma(X_{t})\,dW_{t}
=\int\limits_{0}^{T}\EE\left[\nabla f(X_T)\delta X_T\left|\,X_t\right. \right]\delta X_t^{-1}\sigma(X_t)\,dW_{t},
\end{align}
where $\nabla_xu(t,x)\in\RR^{1\times d}$ denotes the gradient of $u$ w.r.t.~$x$.
The second integral here can be used for
a practical construction of an approximation of $M_T^*$
because the conditional expectation
can be approximated via regression.

The preceding description lacks assumptions
under which the procedure works
(the mentioned ones are not enough).
We refer to~\cite{milstein2009practical} for more detail.

\bigskip
\textbf{Integral representation for the discretisation.}
As was mentioned in the introduction,
we are going to reduce not the variance in $f(X_T)$
but rather the one in $f(X_{\Delta,T})$,
that is, we aim at constructing control variates
directly for the discretised process.
The fine details of the construction
must of course depend on the discretisation scheme.
For the rest of the paper,
we focus on the Euler scheme, that is, we have
\begin{align}
\label{phi_euler}
\Phi_{\Delta}(
       x,y)=x+\mu(x)\Delta+\sigma(x)y.
\end{align}
We define the ``discretised derivative process''
$\delta^i X^k_{t_j,x}(\Delta,t_l)
:=\frac{\partial X^k_{t_j,x}(\Delta,t_l)}{\partial x_i}$,
for $l\ge j$ and $i,k\in\{1,\ldots,d\}$,
where $X^k_{t_j,x}(\Delta,t_l)$
means the $k$-th component of the (Euler) discretisation
for~\eqref{x_sde} started at time~$t_j$ in~$x$
and evaluated at time~$t_l$ ($\ge t_j$),
and use $\delta^i X^k_{\Delta,t_l}$
as an abbreviation of
$\delta^i X^k_{0,x_0}(\Delta,t_l)$.
Assuming $\mu,\sigma\in C^1$,
we get that the process
$(\delta^i X_{\Delta,j\Delta})_{j=1,\ldots,J}$
has the dynamics
\begin{align}
\label{deriv_discr_euler}
\delta^i X_{\Delta, j\Delta} 
    = \delta^i X_{\Delta,(j-1)\Delta}+\sum_{k=1}^d\delta^i X^k_{\Delta,(j-1)\Delta}\left[\frac{\partial\mu(X_{\Delta,(j-1)\Delta})}{\partial x_k}\Delta
+\frac{\partial\sigma(X_{\Delta,(j-1)\Delta})}{\partial x_k}\Delta_jW\right],
\end{align}
(cf.~\eqref{dx_sde}),
where $\delta X_{\Delta,0}=I_d$,
and in what follows $I_d$ denotes the identity matrix of size~$d$.

Given a Borel function $f\colon\bbR^d\to\bbR$
satisfying $\EE\left|f(X_{\Delta,T})\right|<\infty$,
it can be verified by a direct calculation that, for $t\in[t_{j-1},t_j)$,
\begin{align}
\EE[f(X_{\Delta,T}) | \cF_t]=u_\Delta
(t,X_{\Delta,t_{j-1}},W_t-W_{t_{j-1}}),
\label{eq:17022017a1}
\end{align}
where the function $u_{\Delta}\colon\left[0,T\right]\times\bbR^{d+m}\to\bbR$ is constructed via the backward recursion as follows
\begin{align}
u_{\Delta}(t,x,y)&=
\EE[u_{\Delta}(t_j,\Phi_\Delta(x,y+z_j\sqrt{t_j-t}),0)],\quad t\in[t_{j-1},t_j),
\label{eq:17022017a2}\\
u_\Delta(T,x,0)&= f(x),
\label{eq:17022017a3}
\end{align}
where $t_j\doteq\frac{jT}{J}$, $j\in\left\{0,\ldots,J\right\}$, and $z_1,\ldots,z_J\stackrel{\textnormal{i.i.d.}}{\sim}\mathcal{N}(0_m,I_m)$. 

We now introduce the following assumptions:
for any $j\in\{1,\ldots,J\}$ and $x\in\bbR^d$, it holds
\begin{description}
\item[(Ass1)]
$f(X_{t_{j-1},x}(\Delta,T))\in L^1$,
\item[(Ass2)$_n$]
$|\Delta_j W|^n\,\EE[f(X_{t_{j-1},x}(\Delta,T))|\cF_{t_j}]
\in L^1$.
\end{description}
(Ass1)~is just a minimal assumption
that allows to have~\eqref{eq:17022017a1}
with the function $u_\Delta$
constructed via
\eqref{eq:17022017a2}--\eqref{eq:17022017a3}.
(Ass2)$_n$~is a technical assumption,
which depends on~$n$,
allowing to replace integration and differentiation
in several cases of interest (see below).
In most places we need the variant~(Ass2)$_1$,
i.e.\ with $n=1$,
but at a couple of instances we will need
stronger variants (Ass2)$_n$ with $n\ge1$.
That is why we have the parameter~$n$
in the formulation of that assumption.

An attractive feature of such an approach
via the discretised process
(in contrast to the one via the exact solution)
is that, under (Ass1) and~(Ass2)$_1$,
due to the smoothness of the Gaussian density,
the function $u_\Delta$ is continuously differentiable in~$y$
regardless of whether $f$ is smooth,
and, moreover, $u_\Delta$ is continuously differentiable in~$x$,
provided $f,\mu,\sigma$ are continuously differentiable.
More precisely, we obtain the above statements
because, for $t\in[t_{j-1},t_j)$, we can write
(for simplicity, in the one-dimensional case)
$$
u_\Delta(t,x,y)=\int_\bbR
u_\Delta\left(t_j,\Phi_\Delta(x,w),0\right)
\frac1{\sqrt{2\pi(t_j-t)}}
e^{-\frac{(w-y)^2}{2(t_j-t)}}\,dw,
$$
and differentiation under the integral applies
due to~(Ass2)$_1$ together with
the dominated convergence theorem
(notice that the expression
$\EE[f(X_{t_{j-1},x}(\Delta,T))|\cF_{t_j}]$
in~(Ass2)$_n$ is nothing else than
$u_\Delta(t_j,\Phi_\Delta(x,\Delta_j W),0)$).

\begin{theorem}\label{thm:IntegRepr}
Suppose (Ass1) and~(Ass2)$_1$.

\smallskip
(i) It holds
$$
f(X_{\Delta,T})=\EE[f(X_{\Delta,T})]+  
\sum_{j=1}^{J}
\,\int\limits_{t_{j-1}}^{t_j}\nabla_y u_{\Delta}(t,X_{\Delta,t_{j-1}},W_t-W_{t_{j-1}})\,dW_{t},
$$
where $\nabla_yu_\Delta(t,x,y)\in\RR^{1\times m}$ denotes the gradient of $u_\Delta$ w.r.t. $y$. 

\smallskip
(ii) Assume additionally that $f,\mu,\sigma\in C^1$.
Then we also have the alternative representation
$$
f(X_{\Delta,T})=\EE[f(X_{\Delta,T})]+  
\sum_{j=1}^{J}
\,\int\limits_{t_{j-1}}^{t_j}\EE\left[\nabla f(X_{\Delta,T})
\delta X_{\Delta,T}\delta X_{\Delta,t_j}^{-1}\left|\,\mathcal{F}_t\right. \right]\sigma(X_{\Delta,t_{j-1}})\,dW_{t}.
$$
\end{theorem}

Let us define the function $g_j\colon\bbR^d\to\bbR^{1\times d}$, $j\in\left\{1,\ldots,J\right\}$, through
\begin{align}
\label{def_g_strong}
g_j(x)=\left(g_{j,1}(x),\ldots,g_{j,d}(x)\right)\doteq \EE\left[\nabla f(X_{\Delta,T})
\delta X_{\Delta,T}\delta X_{\Delta,t_j}^{-1}\left|\,X_{\Delta,t_{j-1}}=x\right. \right].
\end{align}
Note that it holds (see the proof of Theorem~\ref{thm:IntegRepr})
\begin{align}
&g_{j}(x)=\EE\left[\nabla_x u_\Delta(t_j,X_{\Delta,t_j},0)\left|\,X_{\Delta,t_{j-1}}=x\right. \right],
\label{eq:13022017a1}\\
&\nabla_y u_\Delta(t_{j-1},x,0)=g_j(x)\sigma(x),
\label{eq:13022017a2}
\end{align}
where $\nabla_xu_\Delta(t,x,y)$ denotes the gradient of $u_\Delta$ w.r.t. $x$, and we conditioned
on $X_{\Delta,t_{j-1}}$ instead of $\mathcal{F}_{t_{j-1}}$
because $(X_{\Delta,t_j})_{j=0,\ldots,J}$ is a Markov chain
(one can do that for grid points only).
Theorem~\ref{thm:IntegRepr} inspires to introduce the control variate
\begin{align}
\notag
M_{\Delta,T}^{int,1}\doteq&
\sum_{j=1}^J \sum_{i=1}^m\frac{\partial u_{\Delta}(t_{j-1},X_{\Delta,t_{j-1}},0)}{\partial y_i}\Delta_jW^i\\
\label{cv_deriv}
=&\sum_{j=1}^J\sum_{k=1}^dg_{j,k}(X_{\Delta,t_{j-1}})\sum_{i=1}^m\sigma_{ki}(X_{\Delta,t_{j-1}})\Delta_j W^i.
\end{align}
It will turn out that $M_{\Delta,T}^{int,1}=M_{\Delta,T}^{ser,1}$.
To this end, we derive a connection between the series and integral representations.

\begin{theorem}
\label{relation_eul}
Under (Ass1) and (Ass2)$_n$ for all $n\in\bbN$,
provided that it holds
\begin{align}
\label{deriv_bound}
\left|D^\alpha\left(\frac{\partial}{\partial y_r}u_\Delta(t,x,y)\right)\right|\doteq
\left|\frac{\partial^{K}\left(\frac{\partial}{\partial y_r}u_\Delta(t,x,y)\right)}{\partial t^{\alpha_1}\partial y_1^{\alpha_2}\cdots\partial y_m^{\alpha_{m+1}}}\right|\le C^{K}
\end{align} 
for all $K\in\NN$, $r\in\left\{1,\ldots,m\right\}$, $\left|\alpha\right|=K$, $t\in[t_{j-1},t_j)$, $x\in\RR^d$, $y\in\RR^m$,
with some constant $C>0$, we have for the Euler scheme
\begin{align}
\label{eul_cv_ser}
f(X_{\Delta,T})=\EE[f(X_{\Delta,T})]+ \sum_{j=1}^J\sum_{l=1}^\infty\Delta^{l/2}\sum_{\substack{k\in\NN_0^m\\ \sum_{r=1}^mk_r=l}}\frac{\partial^lu_{\Delta}(t_{j-1},X_{\Delta,t_{j-1}},0)}{\partial y_1^{k_1}\cdots\partial y_m^{k_m}}\prod_{r=1}^m\frac{H_{k_r}\left(\frac{\Delta_jW^r}{\sqrt{\Delta}}\right)}{\sqrt{k_r!}}
\end{align} 
whenever $0<\Delta<\frac{1}{C^2}$.
(The series converge in~$L^2$.)
Consequently, we obtain for $l=\sum_{r=1}^m k_r\in\NN$
\begin{align}
\label{cv_equiv}
\frac{\Delta^{l/2}}{\sqrt{k_1!}\cdots \sqrt{k_m!}}\cdot\frac{\partial^lu_{\Delta}(t_{j-1},X_{\Delta,t_{j-1}},0)}{\partial y_1^{k_1}\cdots\partial y_m^{k_m}}=a_{j,k}(X_{\Delta, t_{j-1}}).
\end{align}
\end{theorem}

\begin{remark}
In the one-dimensional case ($d=m=1$), a representation of a similar type as~\eqref{eul_cv_ser}
appears in~\cite{AAO} in a somewhat different form.
Our form is aimed
at constructing control variates
via regression methods.
\end{remark}

In particular, we see from Theorem~\ref{relation_eul}
that $M_{\Delta,T}^{int,1}=M_{\Delta,T}^{ser,1}$
provided that~\eqref{deriv_bound} holds.
However, we can prove the equality of
the aforementioned control variates
without assuming~\eqref{deriv_bound}:

\begin{theorem}
\label{th:13022017a1}
Under (Ass1) and (Ass2)$_1$,
we have for $i\in\left\{1,\ldots,m\right\}$
\begin{align*}
a_{j,e_i}(x)=\sqrt{\Delta}\frac{\partial}{\partial y_i}u_\Delta(t_{j-1},x,0),
\end{align*}
and consequently,
\begin{align*}
M_{\Delta,T}^{int,1}=M_{\Delta,T}^{ser,1}.
\end{align*}
\end{theorem}

It is interesting to remark that,
although we assumed $f(X_{\Delta,T})\in L^2$
when speaking about the series representation,
the coefficients $a_{j,e_i}$ are well-defined
already under (Ass1) and~(Ass2)$_1$.

We can now investigate the order of the truncation error,
which arises when we
replace the control variate $M_{\Delta,T}$ of~\eqref{cv_perfect}
with the control variate $M_{\Delta,T}^{ser,1}$ of~\eqref{cv_truncated}.

\begin{theorem}
\label{var_ord_eul}
Suppose (Ass1) and~(Ass2)$_3$.
Provided that the function $u_\Delta(t,x,y)$
has bounded partial derivatives in $y$ of orders 2 and 3, it holds
\begin{align}
\label{cv_var_ord}
\Var\left[f(X_{\Delta, T})-M_{\Delta,T}^{int,1}\right]=\Var\left[f(X_{\Delta, T})-M_{\Delta,T}^{ser,1}\right]\lesssim\Delta.
\end{align}
\end{theorem}

\begin{remark}
(i) Below we will present sufficient conditions in terms of
the functions $f,\mu,\sigma$
that ensure the assumption on $u_\Delta$ in Theorem~\ref{var_ord_eul}
(see Theorem~\ref{func_assump} in Section~\ref{sec:3}).

\smallskip
(ii) The control variate $M_{\Delta,T}^{int,1}$ 
differs from the one suggested in~\cite{milstein2009practical} only in an index concerning the inverted matrix, i.e. we have $\delta X_{\Delta,t_j}^{-1}$ inside of $g_j(X_{\Delta,t_{j-1}})$ rather than the $\mathcal{F}_{t_{j-1}}$-measurable random variable $\delta X_{\Delta,t_{j-1}}^{-1}$ which arises in case of the exact solution $f(X_T)$ from a simple discretisation of the stochastic integral
in~\eqref{cv_milstein}.
\end{remark}

Regarding the weak convergence order of the Euler scheme, we have the following result
(cf.~Theorem~2.1 in~\cite{MilsteinTretyakov:2004}).

\begin{proposition}
\label{prop:17022017a1}
Assume that $\mu$ and $\sigma$ in~\eqref{x_sde} are Lipschitz continuous with components $\mu_k,\,\sigma_{ki}\colon \RR^d\to\RR$,
$k=1,\ldots,d$, $i=1,\ldots,m$,
being $4$ times continuously differentiable
with their partial derivatives of orders up to $4$
having polynomial growth.
Let $f\colon\RR^d\to\RR$ be $4$ times continuously differentiable with 
partial derivatives of orders up to $4$
having polynomial growth.
Then, for the Euler scheme~\eqref{phi_euler}, we have
\begin{align}
\label{bias_euler}
\left|\EE f(X_T)-\EE f(X_{\Delta,T})\right|\le c\Delta,
\end{align}
where the constant $c$ does not depend on $\Delta$. 
\end{proposition}

We remark that the assumption that,
for sufficiently large $n\in\mathbb N$,
the expectations $\EE |X_{\Delta,j\Delta}|^{2n}$
are uniformly bounded in $J$ and $j=0,\ldots,J$
(cf.~Theorem~2.1 in~\cite{MilsteinTretyakov:2004})
is automatically satisfied for the Euler scheme
because $\mu$ and $\sigma$, being globally Lipschitz,
have at most linear growth.

\section{Regression analysis}
\label{sec:3}
In the previous sections we have given several representations for the control variates. Now  we discuss how to compute the coefficients in these representations via regression. For the sake of clarity, we will focus on the control variate given by~\eqref{cv_deriv}, that is, we will estimate the functions $g_{j,k}$ in~\eqref{def_g_strong} via linear regression.
Let us start with a general description of the global Monte Carlo regression algorithm.

\subsection{Global Monte Carlo regression algorithm}
\label{glob_mc_reg_strong}
Fix a $q$-dimensional vector of real-valued functions \ensuremath{\psi=(\psi^{1},\ldots,\psi^{q})}
 on \ensuremath{\RR^{d}}. Simulate a set of $N$ ``training paths'' of the Markov chains $X_{\Delta,j\Delta}$ and $\delta X_{\Delta,j\Delta}$, $j=0,\ldots, J$.
We should choose $N>q$. In what follows these $N$ training paths
are denoted by $D_N^{tr}$:
$$
D_N^{tr}\doteq
\left\{
(X_{\Delta,j\Delta}^{tr,(n)},\delta X_{\Delta,j\Delta}^{tr,(n)})_{j=0,\ldots,J}:
n=1,\ldots,N
\right\}.
$$
Let $\boldsymbol{\alpha}_{j,k}=(
\alpha_{j,k}^{1},\ldots,\alpha_{j,k}^{q})$,
where $j\in\left\{1,\ldots, J\right\}$, $k\in\left\{1,\ldots,d\right\}$,
 be a solution of the following least squares optimisation problem:
\begin{align*}
\operatorname{argmin}_{\boldsymbol{\alpha}\in\RR^{q}}
\sum_{n=1}^{N}\left[\zeta^{tr,(n)}_{j,k}-\alpha^{1}\psi^{1}(X_{\Delta,(j-1)\Delta}^{tr,(n)})-\ldots-\alpha^{q}\psi^{q}(X_{\Delta,(j-1)\Delta}^{tr,(n)})\right]^{2}
\end{align*}
with 
\begin{align*}
\zeta_{j}^{tr,(n)}=\left(\zeta_{j,1}^{tr,(n)},\ldots,\zeta_{j,d}^{tr,(n)}\right)\doteq \nabla f(X_{\Delta,T}^{tr,(n)})
\delta X_{\Delta,T}^{tr,(n)}\left(\delta X_{\Delta,j\Delta}^{tr,(n)}\right)^{-1}.  
\end{align*}
Define an estimate for  the coefficient function $g_{j,k}$ via
\begin{align*}
\hat g_{j,k}(z)\doteq \alpha_{j,k}^{1}\psi^{1}(z)+\ldots+\alpha_{j,k}^{q}\psi^{q}(z),\quad z\in\RR^{d}.
\end{align*}
The cost of computing
$\boldsymbol{\alpha}_{j,k}$ is of order $O(Nq^{2})$,
 since each \ensuremath{\boldsymbol{\alpha}_{j,k}}
 is of the form \ensuremath{\boldsymbol{\alpha}_{j,k}=B^{-1}b}
 with 
\begin{align}
\label{strong_b_matr_reg}
B_{l,o}\doteq\frac{1}{N}\sum_{n=1}^{N}\psi^{l}\bigl(X_{\Delta,(j-1)\Delta}^{tr,(n)}\bigr)\psi^{o}\bigl(X_{\Delta,(j-1)\Delta}^{tr,(n)}\bigr)
\end{align}
and 
\begin{align*}
b_{l}\doteq\frac{1}{N}\sum_{n=1}^{N}\psi^{l}\bigl(X_{\Delta,(j-1)\Delta}^{tr,(n)}\bigr)\,\zeta_{j,k}^{tr,(n)},
\end{align*}
\ensuremath{l,o\in\{1,\ldots,q\}.} The cost of approximating the family of the coefficient functions $g_{j,k}$, $j\in\left\{1,\ldots, J\right\}$, $k\in\left\{1,\ldots,d\right\}$, is of order
$O\bigl(JdNq^{2}\bigr)$.

\subsection{Piecewise polynomial regression}
\label{piece_poly}
There are different ways to choose the basis functions \ensuremath{\psi=(\psi^{1},\ldots,\psi^{q}).}
In this section we describe
piecewise polynomial partitioning estimates
and present $L^2$-upper bounds
for the estimation error.

From now on, we fix some $p\in\bbN_0$,
which will denote the maximal degree
of polynomials involved in our basis functions.
The piecewise polynomial partitioning estimate of $g_{j,k}$ works as follows: 
consider some $R>0$ and an equidistant partition of $\left[-R,R\right]^d$ in $Q^d$ cubes $K_1,\ldots,K_{Q^d}$. Further, consider the basis functions $\psi^{l,1},\ldots,\psi^{l,q}$ with $l\in\left\{1,\ldots,Q^d\right\}$ and $q=\binom{p+d}{d}$ such that $\psi^{l,1}(x),\ldots,\psi^{l,q}(x)$ are polynomials with degree less than or equal to $p$ for $x\in K_l$ and $\psi^{l,1}(x)=\ldots=\psi^{l,q}(x)=0$ for $x\notin K_l$. Then we obtain the least squares regression estimate $\hat g_{j,k}(x)$ for $x\in\RR^d$ as described in Section~\ref{glob_mc_reg_strong}, based on $Q^dq=O(Q^dp^d)$ basis functions. In particular, we have $\hat g_{j,k}(x)=0\) for any \(x\notin\left[-R,R\right]^d$. We note that the cost of computing $\hat g_{j,k}$ for all $j,k$ is of order $O(Jd NQ^{d} p^{2d})$ rather than $O(Jd N Q^{2d} p^{2d})$ due to a block diagonal matrix structure
of $B$ in~\eqref{strong_b_matr_reg}. An equivalent approach, which leads to the same estimator $\hat g_{j,k}(x)$, is to perform separate regressions for each cube $K_1,\ldots,K_{Q^d}$. Here, the number of basis functions at each regression is of order $O(p^d)$ so that the overall cost is of order $O(Jd N Q^{d} p^{2d})$, too.
For $x=(x_1,\ldots,x_d)\in\bbR^d$ and $h\in[1,\infty)$,
we will use the notations
\begin{align*}
|x|_h\doteq\bigg(\sum_{k=1}^d |x_k|^h\bigg)^{1/h},\quad
|x|_\infty\doteq\max_{k=1,\ldots,d}|x_k|.
\end{align*}
For $s\in\bbN_0$,
$C>0$ and $h\in[1,\infty]$,
we say that
a function $F\colon\bbR^d\to\bbR$ is
\emph{${(s+1,C)}$-smooth w.r.t.\ the norm $\left|\cdot\right|_h$} whenever, for all
$\alpha=(\alpha_1,\ldots,\alpha_d)\in\bbN_0^d$
with $\sum_{k=1}^d \alpha_k=s$, we have
\begin{align*}
|D^\alpha F(x)-D^\alpha F(y)|\le C|x-y|_h,
\quad x,y\in\bbR^d,
\end{align*}
i.e.\ the function $D^\alpha F$
is globally Lipschitz
with the Lipschitz constant $C$
with respect to the norm $|\cdot|_h$
on $\bbR^d$
(cf.~Definition~3.3 in~\cite{gyorfi2002distribution}).
In what follows,
we use the notation
$\PP_{\Delta,j-1}$
for the distribution of $X_{\Delta,(j-1)\Delta}$. In particular, we will work with the corresponding $L^2$-norm:
\begin{align*}
\|F\|^2_{L^2(\PP_{\Delta,j-1})}\doteq \int_{\RR^d} F^2(x)\,\PP_{\Delta,j-1}(dx)=\mathbb{E}\left[F^2\left(X_{\Delta,(j-1)\Delta}\right)\right].
\end{align*}
We now define $\zeta_{j,k}$ as the $k$-th component of the vector $\zeta_j=\left(\zeta_{j,1},\ldots,\zeta_{j,d}\right)\doteq \nabla f(X_{\Delta,T})\delta X_{\Delta,T}\delta X_{\Delta,j\Delta}^{-1}$
and remark that
$g_{j,k}(x)=\EE[\zeta_{j,k}|X_{\Delta, (j-1)\Delta}=x]$.
In what follows, we consider the following assumptions:
there exist $h\in[1,\infty]$
and positive constants $\Sigma,A,C_h,\nu,B_\nu$
such that, for all $J\in\bbN$,
$j\in\left\{1,\ldots, J\right\}$ and $k\in\left\{1,\ldots,d\right\}$,
it holds
\begin{itemize}
\item[(A1)]
$\quad\sup_{x\in\RR^d}\Var[\zeta_{j,k}|X_{\Delta,(j-1)\Delta}=x]
\le\Sigma<\infty$,
\item[(A2)]
$\quad\sup_{x\in\RR^d} |g_{j,k}(x)|\le A<\infty$,
\item[(A3)]
$\quad g_{j,k}$ is $(p+1,C_h)$-smooth w.r.t.\ the norm $|\cdot|_h$,
\item[(A4)]
$\quad\PP(|X_{\Delta,(j-1)\Delta}|_\infty>R)\le B_\nu R^{-\nu}$ for all $R>0$.
\end{itemize}

\begin{remark}
Let us notice that it is only a matter of convenience
which $h$ to choose in~(A3)
because all norms $|\cdot|_h$ are equivalent.
Furthermore, since $\mu$ and $\sigma$ are assumed
to be globally Lipschitz, hence have linear growth,
then, given any $\nu>0$,
(A4)~is satisfied with a sufficiently large $B_\nu>0$.
In other words, (A4)~is needed
only to introduce the constant $B_\nu$,
which appears in the formulations below.
\end{remark}

In the next theorem we, in particular, present sufficient conditions
in terms of the functions $\mu$, $\sigma$ and~$f$
that imply the preceding assumptions.

\begin{theorem}
\label{func_assump}
(i) Under (Ass1) and (Ass2)$_1$,
let all functions $f,\mu_k,\sigma_{ki}$, $k\in\left\{1,\ldots,d\right\}$, $i\in\left\{1,\ldots,m\right\}$, be continuously differentiable with bounded partial derivatives.
Then (A1) and~(A2) hold with appropriate constants
$\Sigma$ and~$A$.

\smallskip
(ii) If, moreover, (Ass1) and (Ass2)$_3$ are satisfied,
all functions $\sigma_{ki}$ are bounded
and all functions $f,\mu_k,\sigma_{ki}$
are $3$ times continuously differentiable with
bounded partial derivatives up to order~$3$,
then the function $u_\Delta(t,x,y)$
has bounded partial derivatives in $y$ up to order~$3$.
In particular, \eqref{cv_var_ord}~holds true.
\end{theorem}

\begin{remark}
As a generalisation of Theorem~\ref{func_assump}, it is natural to expect that (A3) is satisfied
with a sufficiently large constant $C_h>0$
if, under (Ass1) and (Ass2)$_{p+2}$,
all functions $f,\mu_k,\sigma_{ki}$
are $p+2$ times continuously differentiable
with bounded partial derivatives up to order $p+2$.
\end{remark}

Let $\hat g_{j,k}$ be the piecewise
polynomial partitioning estimate
of~$g_{j,k}$.
By $\tilde g_{j,k}$ we denote the truncated estimate,
which is defined as
\begin{align*}
\tilde g_{j,k}(x)\doteq
T_{A}\hat g_{j,k}(x)
\doteq\begin{cases}
\hat g_{j,k}(x)&\text{if }
|\hat g_{j,k}(x)|\le A,\\
A\sgn\hat g_{j,k}(x)
&\text{otherwise,}
\end{cases}
\end{align*}
where $A$ is the bound from~(A2).

\begin{lemma}
\label{lemm_cv_strong}
Under (A1)--(A4), we have
\begin{align}
\label{norm_lemma_strong}
\EE\|\tilde g_{j,k}-g_{j,k}\|^2_{L^2(\PP_{\Delta,j-1})}
&\le
\tilde c\left(\Sigma+A^2(\log N+1)\right)\frac{\binom{p+d}d Q^d}{N}
\\
\notag
&\hspace{1em}+\frac{8\,C_h^2}{(p+1)!^2 d^{2-2/h}}
\left(\frac{Rd}Q\right)^{2p+2}
+8A^2B_\nu R^{-\nu},
\end{align}
where $\tilde c$ is a universal constant.
\end{lemma}

It is worth noting that the expectation
in the left-hand side of~\eqref{norm_lemma_strong}
accounts for the averaging over the randomness
in $D_N^{tr}$. To explain this in more detail,
let $(X_{\Delta,j\Delta})_{j=0,\ldots,J}$
be a ``testing path'' which is independent
of the training paths $D_N^{tr}$. Then it holds
\begin{align*}
\|\tilde g_{j,k}-g_{j,k}\|^2_{L^2(\PP_{\Delta,j-1})}
&\equiv
\|\tilde g_{j,k}(\cdot,D_N^{tr})-g_{j,k}(\cdot)\|^2_{L^2(\PP_{\Delta,j-1})}\\
&=
\EE\left[
\left(\tilde g_{j,k}(X_{\Delta,(j-1)\Delta},D_N^{tr})
-g_{j,k}(X_{\Delta,(j-1)\Delta})\right)^2
\,|\,D_N^{tr}\right],
\end{align*}
hence,
\begin{align}\label{eq:30042016a4}
\EE\|\tilde g_{j,k}-g_{j,k}\|^2_{L^2(\PP_{\Delta,j-1})}
=\EE\left[
\left(\tilde g_{j,k}(X_{\Delta,(j-1)\Delta},D_N^{tr})
-g_{j,k}(X_{\Delta,(j-1)\Delta})\right)^2
\right],
\end{align}
which provides an alternative form for
the expression in the left-hand side
of~\eqref{norm_lemma_strong}.
\par
Let us now estimate the variance of the random variable
\(f(X_{\Delta,T})-\tilde{M}^{{int},1}_{\Delta,T}\),
where
\begin{align}
\label{eq:cv_strong2}
\tilde{M}^{{int},1}_{\Delta,T}\doteq\sum_{j=1}^J\sum_{k=1}^d\tilde{g}_{j,k}(X_{\Delta,(j-1)\Delta},D_N^{tr})\sum_{i=1}^m\sigma_{ki}(X_{\Delta,(j-1)\Delta})\Delta_j W^i.
\end{align}

\begin{theorem}
\label{theo_cv_strong}
Let us assume $\sup_{x\in\RR^d} |\sigma_{ki}(x)|\le \sigma_{\max}<\infty$
for all $k\in\left\{1,\ldots,d\right\}$ and $i\in\left\{1,\ldots,m\right\}$.
Then we have under (A1)--(A4)
\begin{eqnarray}
\notag
\Var[f(X_{\Delta,T})-
\tilde{M}^{{int},1}_{\Delta,T}]&\lesssim&
\frac{1}{J}+d^2Tm
\sigma_{\max} ^2
\left\{\tilde c\left(\Sigma+A^2(\log N+1)\right)\frac{\binom{p+d}d Q^d}{N}\right.
\\
\label{eq:cv_bound_strong}
&&
\left.+\frac{8\,C_h^2}{(p+1)!^2 d^{2-2/h}}
\left(\frac{Rd}Q\right)^{2p+2}
+8A^2 B_\nu R^{-\nu}\right\}.
\end{eqnarray}
\end{theorem}

We finally stress that
$\tilde{M}^{{int},1}_{\Delta,T}$
is a valid control variate
in that is does not introduce bias,
i.e.\ $\EE[\tilde{M}^{{int},1}_{\Delta,T}|D_N^{tr}]=0$,
which follows from
the martingale transform structure in~\eqref{eq:cv_strong2}.

\subsection{Summary of the algorithm}
The algorithm of the ``integral approach'' consists of two phases:
training phase and testing phase.
In the training phase, we simulate
$N$ independent training paths $D_N^{tr}$
and construct regression estimates
$\tilde g_{j,k}(\cdot,D_N^{tr})$
for the coefficients $g_{j,k}(\cdot)$, $k\in\left\{1,\ldots,d\right\}$.
In the testing phase,
independently from $D_N^{tr}$
we simulate $N_0$ independent testing paths
$(X_{\Delta,j\Delta}^{(n)})_{j=0,\ldots,J}$,
$n=1,\ldots,N_0$,
and build the Monte Carlo estimator
for $\EE f(X_T)$ as
\begin{equation}
\label{eq:cv_testing_strong}
\frac1{N_0}
\sum_{n=1}^{N_0}
\left(f(X^{(n)}_{\Delta,T})-\tilde{M}^{{int},1,(n)}_{\Delta,T}\right).
\end{equation}
The expectation of this estimator
equals $\EE f(X_{\Delta,T})$,
and the upper bound for the variance
is $\frac1{N_0}$ times the expression
in~\eqref{eq:cv_bound_strong}.

\section{Complexity analysis}
\label{sec:4}
The results presented in previous sections
provide us with ``building blocks''
to perform the complexity analysis.

\medskip\noindent
\textbf{Standing assumption for Complexity Analysis}
consists in
\begin{center}
(Ass1), (Ass2)$_1$, \eqref{cv_var_ord} and~\eqref{bias_euler}.
\end{center}
Combining Theorem~\ref{var_ord_eul},
Theorem~\ref{func_assump} and
Proposition~\ref{prop:17022017a1},
we recall that this standing assumption is satisfied
whenever we have (Ass1), (Ass2)$_3$,
$\sigma$ is bounded,
$f,\mu,\sigma\in C^4$,
the partial derivatives
of $f$, $\mu$ and $\sigma$
up to order~3 are bounded
and of order~4 have polynomial growth.
However, we prefer to formulate the standing assumption
for complexity analysis as above because
one might imagine other sufficient conditions for~it.

\subsection{Integral approach}
\label{sec:compl_strong}
Below we present a complexity analysis
which explains how we can approach
the complexity order $\varepsilon^{-2}\sqrt{\left|\log(\varepsilon)\right|}$
with $\varepsilon$ being the precision
to be achieved.

For the integral approach we perform $d$ regressions in the training phase
and $d$ evaluations of $\tilde g_{j,k}$
in the testing phase (using the regression coefficients
from the training phase) at each time step.
Therefore, the overall cost is of order
\begin{align}
\label{cost_strong}
JQ^ddc_{p,d}\max\left\{c_{p,d}N,N_0\right\},
\end{align}
where $c_{p,d}\doteq\binom{p+d}{p}$.
Under (A1)--(A4) and boundedness of~$\sigma$
(cf.~Theorem~\ref{theo_cv_strong}),
we have the following constraints
\begin{align}
\label{mse_strong_constr}
\max\left\{\frac{1}{J^2},\frac{1}{JN_0},
\frac{Q^dd^2mc_{p,d}\log(N)}{NN_0},
\frac{d^2m}{(p+1)!^2N_0}\left(\frac{Rd}{Q}\right)^{2(p+1)},
\frac{d^2mB_\nu}{N_0R^\nu}\right\}\lesssim\varepsilon^2,
\end{align} 
to ensure a \textit{mean squared error} (MSE) of order $\varepsilon^2$.
Note that the first term in~\eqref{mse_strong_constr} comes from the squared bias of the estimator (due to~\eqref{bias_euler} and $\EE[\tilde M^{{int},1}_{\Delta,T}]=0$) and the remaining four ones come from the variance of the estimator (see~\eqref{eq:cv_bound_strong} and~\eqref{eq:cv_testing_strong}). 

\begin{theorem}
\label{theo_comp_strong}
Under (A1)--(A4) and boundedness of~$\sigma$,
we obtain the following solution for the integral approach:
\begin{align}
\label{J_strong}
J&\asymp \varepsilon^{-1},\quad 
Q\asymp \left[\frac{B_\nu^{4(p+1)}d^{2\nu+4(p+1)(\nu+1)}m^{\nu+2(p+1)}}{\varepsilon^{2\nu+4(p+1)}c_{p,d}^{2\nu+4(p+1)}(p+1)!^{4\nu}}\right]^\frac{1}{d\nu+2(p+1)(d+2\nu)}, \\
\notag
N&\asymp\left[\frac{B_\nu^{2d(p+1)}d^{2d\nu+2(p+1)(d\nu+2d+2\nu)}m^{d\nu+2(p+1)(d+\nu)}}{\varepsilon^{2d\nu+4(p+1)(d+\nu)}c_{p,d}^{d\nu+2d(p+1)}(p+1)!^{2d\nu}}\right]^\frac{1}{d\nu+2(p+1)(d+2\nu)}\\
\label{N_strong}
&\phantom{\asymp}\cdot\sqrt{\log\left(\varepsilon^{-\frac{2d\nu+4(p+1)(d+\nu)}{d\nu+2(p+1)(d+2\nu)}}\right)}, \\
\notag
 N_0&\asymp Nc_{p,d}\\
 \notag
 &\asymp \left[\frac{B_\nu^{2d(p+1)}c_{p,d}^{4\nu(p+1)}d^{2d\nu+2(p+1)(d\nu+2d+2\nu)}m^{d\nu+2(p+1)(d+\nu)}}{\varepsilon^{2d\nu+4(p+1)(d+\nu)}(p+1)!^{2d\nu}}\right]^\frac{1}{d\nu+2(p+1)(d+2\nu)}\\
  \label{N0_strong}
 &\phantom{\asymp}\cdot\sqrt{\log\left(\varepsilon^{-\frac{2d\nu+4(p+1)(d+\nu)}{d\nu+2(p+1)(d+2\nu)}}\right)}, \\
 R&\asymp \left[\frac{B_\nu^{d+4(p+1)}(p+1)!^{2d}m^{2(p+1)}}{\varepsilon^{4(p+1)}c_{p,d}^{4(p+1)}d^{2(p+1)(d-2)}}\right]^\frac{1}{d\nu+2(p+1)(d+2\nu)}, 
\end{align}
provided that $2(p+1)>d$ and $\nu>\frac{2d(p+1)}{2(p+1)-d}$.\footnote{\label{ft:compl_strong}Performing 
the full complexity analysis via Lagrange multipliers
one can see that these parameter values are
\emph{not} optimal if $2(p+1)\le d$ or $\nu\le\frac{2d(p+1)}{2(p+1)-d}$
(a Lagrange multiplier corresponding to
a ``$\le0$'' constraint is negative, cf. proof of Theorem~\ref{theo_comp_strong}).
Therefore, the recommendation is to choose
the power $p$ for our basis functions
according to $p>\frac{d-2}2$.
The opposite choice is allowed as well
(the method converges),
but theoretical complexity of the method
would be then worse than that of the SMC,
namely, $\varepsilon^{-3}$.}
Thus, we have for the complexity
\begin{align}
\notag
\mathcal{C}_{int}&\asymp JQ^ddc_{p,d}^2N\asymp JQ^ddc_{p,d}N_0\\
\notag
&\asymp\left[\frac{B_\nu^{6d(p+1)}c_{p,d}^{2(p+1)(4\nu-d)-d\nu}d^{5d\nu+2(p+1)(3d\nu+5d+4\nu)}m^{3d\nu+6(p+1)(d+\nu)}}{\varepsilon^{5d\nu+2(p+1)(5d+4\nu)}(p+1)!^{6d\nu}}\right]^\frac{1}{d\nu+2(p+1)(d+2\nu)}\\
\label{compl_strong}
&\phantom{\asymp}\cdot\sqrt{\log\left(\varepsilon^{-\frac{2d\nu+4(p+1)(d+\nu)}{d\nu+2(p+1)(d+2\nu)}}\right)}.
\end{align}
\end{theorem}

\begin{remark}
(i) For the sake of comparison with the SMC and MLMC approaches,
we recall at this point that their complexities are
$$
\mathcal C_{SMC}\asymp\varepsilon^{-3}
\quad\text{and}\quad
\mathcal C_{MLMC}\asymp\varepsilon^{-2}
$$
at best.\footnote{For the Euler scheme,
there is an additional logarithmic factor in the complexity of the MLMC algorithm (see~\cite{giles2008multilevel}).}
Complexity estimate~\eqref{compl_strong}
shows that one can approach the complexity order
$\varepsilon^{-2}\sqrt{\left|\log(\varepsilon)\right|}$, when $p,\nu\to\infty$,
i.e.\ if the coefficients \(g_{j,k}\) are smooth enough and the solution $X$ of SDE~\eqref{x_sde} lives in a compact set.

\smallskip
(ii) Note that we would have obtained the same complexity even when the variance in~\eqref{cv_var_ord} were of order $\Delta^K$ with $K>1$. This is due to the fact that the second constraint in~\eqref{mse_strong_constr} is the only inactive one and this would still hold if the condition were $\frac{1}{J^KN_0}\lesssim\varepsilon^2$. Hence, it is not useful to derive a control variate with a higher variance order for the Euler scheme.
\end{remark}

\subsection{Series approach}
\label{sec:compl_strong_ser}
Below we present a complexity analysis for the series representation, defined in Section~\ref{series}. Again we focus on the Euler scheme~\eqref{phi_euler}. Then we compare the resulting complexity with the one in~\eqref{compl_strong}.

Similarly to Section~\ref{piece_poly},
we define $\zeta_{j,i}$ as the $i$-th component of the vector $\zeta_j=\left(\zeta_{j,1},\ldots,\zeta_{j,m}\right)^\top\doteq f(X_{\Delta,T})\frac{\Delta_j W}{\sqrt{\Delta}}$
and remark that
$a_{j,e_i}(x)=\EE[\zeta_{j,i}|X_{\Delta, (j-1)\Delta}=x
]$ (compare with~\eqref{eq:wiener_pcv}).
We will work under the following assumptions:
there exist $h\in[1,\infty]$
and positive constants $\Sigma,A,C_h$
such that, for all $J\in\bbN$,
$j\in\left\{1,\ldots, J\right\}$ and $i\in\left\{1,\ldots,m\right\}$,
it holds:
\begin{itemize}
\item[(B1)]
$\quad\sup_{x\in\RR^d}\Var[\zeta_{j,i}|X_{\Delta,(j-1)\Delta}=x]
\le\Sigma<\infty$,
\item[(B2)]
$\quad\sup_{x\in\RR^d} |a_{j,e_i}(x)|\le A\sqrt{\Delta}<\infty$,
\item[(B3)]
$\quad a_{j,e_i}$ is $(p+1,C_h)$-smooth w.r.t.\ the norm $|\cdot|_h$.
\end{itemize}
Note the difference between~(B2)
and (A2) of Section~\ref{piece_poly},
while (B1) has the same form as~(A1).
This is due to~\eqref{cv_equiv},
hence the additional factor $\sqrt{\Delta}$ in~(B2).

In what follows the $N$ training paths
are denoted by 
$$
D_N^{tr}\doteq
\left\{
(X_{\Delta,j\Delta}^{tr,(n)})_{j=0,\ldots,J}:
n=1,\ldots,N
\right\},
$$
that is, we do not need to simulate paths for the derivative processes $\delta X_{\Delta,j\Delta}$. 
Let $\hat a_{j,e_i}$ be the piecewise
polynomial partitioning estimate
of $a_{j,e_i}$ described
in Section~\ref{piece_poly}.
By $\tilde a_{j,e_i}$ we denote the truncated estimate,
which is defined as follows: 
\begin{align*}
\tilde a_{j,e_i}(x)\doteq
T_{A\sqrt{\Delta}}\hat a_{j,e_i}(x)
\doteq\begin{cases}
\hat a_{j,e_i}(x)&\text{if }
|\hat a_{j,e_i}(x)|\le A\sqrt{\Delta},\\
A\sqrt{\Delta}\sgn\hat a_{j,e_i}(x)
&\text{otherwise.}
\end{cases}
\end{align*}

\begin{lemma}
\label{lemm_cv_strong_ser}
Under (B1)--(B3) and (A4), we have
\begin{align}
\label{norm_lemma_strong_ser}
\EE\|\tilde a_{j,e_i}-a_{j,e_i}\|^2_{L^2(\PP_{\Delta,j-1})}
&\le
\tilde c\left(\Sigma+A^2\Delta(\log N+1)\right)\frac{c_{p,d} Q^d}{N}
\\
\notag
&\hspace{1em}+\frac{8\,C_h^2}{(p+1)!^2d^{2-\frac{2}{h}}}
\left(\frac{R}Q\right)^{2p+2}
+8A^2\Delta B_\nu R^{-\nu},
\end{align}
where $\tilde c$ is a universal constant.
\end{lemma}
\par
Let us now estimate the variance of the random variable
\(f(X_{\Delta,T})-\tilde{M}^{{ser},1}_{\Delta,T}\),
where 
\begin{align}
\label{eq:cv_strong_ser}
\tilde{M}^{{ser},1}_{\Delta,T}\doteq\sum_{j=1}^J\sum_{i=1}^m\tilde{a}_{j,e_i}(X_{\Delta,(j-1)\Delta},D_N^{tr})\frac{\Delta_j W^i}{\sqrt{\Delta}}.
\end{align}
\begin{theorem}
\label{theo_cv_strong_ser}
Under (B1)--(B3) and (A4), we have
\begin{eqnarray}
\notag
\Var[f(X_{\Delta,T})-
\tilde{M}^{{ser},1}_{\Delta,T}]&\lesssim&
\frac{1}{J}+Jm
\left\{\tilde c\left(\Sigma+A^2\Delta(\log N+1)\right)\frac{c_{p,d}Q^d}{N}\right.
\\
\label{eq:cv_bound_strong_ser}
&&
\left.+\frac{8\,C_h^2}{(p+1)!^2d^{2-\frac{2}{h}}}
\left(\frac{R}Q\right)^{2p+2}
+8A^2\Delta B_\nu R^{-\nu}\right\}.
\end{eqnarray}
\end{theorem}
Let us study the complexity of the following ``series approach'': In the training phase, we simulate
$N$ independent training paths $D_N^{tr}$
and construct regression estimates
$\tilde a_{j,e_i}(\cdot,D_N^{tr})$
for the coefficients $a_{j,e_i}(\cdot)$, $i\in\left\{1,\ldots,m\right\}$.
In the testing phase,
independently from $D_N^{tr}$
we simulate $N_0$ independent testing paths
$(X_{\Delta,j\Delta}^{(n)})_{j=0,\ldots,J}$,
$n=1,\ldots,N_0$,
and build the Monte Carlo estimator
for $\EE f(X_T)$ as
\begin{equation}
\label{eq:cv_testing_strong_ser}
\frac1{N_0}
\sum_{n=1}^{N_0}
\left(f(X^{(n)}_{\Delta,T})-\tilde M^{{ser},1,(n)}_{\Delta,T}\right).
\end{equation}
Therefore, the overall cost is of order
\begin{align}
\label{cost_strong_ser}
JQ^dmc_{p,d}\max\left\{c_{p,d}N,N_0\right\}.
\end{align}
The expectation of the estimator in~\eqref{eq:cv_testing_strong_ser}
equals $\EE f(X_{\Delta,T})$,
and the upper bound for the variance
is $\frac1{N_0}$ times the expression
in~\eqref{eq:cv_bound_strong_ser}.
Hence, we have the following constraints 
\begin{align}
\label{mse_strong_constr_ser}
\max\left\{\frac{1}{J^2},\frac{1}{JN_0},
\frac{JQ^dmc_{p,d}}{NN_0},
\frac{Jm}{(p+1)!^2N_0}\left(\frac{Rd}{Q}\right)^{2(p+1)},
\frac{mB_\nu}{N_0R^\nu}\right\}\lesssim\varepsilon^2,
\end{align} 
to ensure a MSE of order $\varepsilon^2$ (due to $\EE[M_{\Delta,T}^{{ser},1}]=0$ as well as ~\eqref{eq:cv_bound_strong_ser} and~\eqref{eq:cv_testing_strong_ser}). Note that there is no longer a $\log$-term in~\eqref{mse_strong_constr_ser}. This is due to the factor $\Delta$ in~\eqref{eq:cv_bound_strong_ser} such that $\Sigma$ is of a higher order, compared to $\Delta(\log N+1)$.

\begin{theorem}
\label{theo_comp_strong_ser}
Under (B1)--(B3) and~(A4),
we obtain the following solution for the series approach:
\begin{align}
\label{J_strong_ser}
J&\asymp \varepsilon^{-1},\quad 
Q\asymp \left[\frac{B_\nu^{4(p+1)}d^{4\nu(p+1)}m^{\nu+2(p+1)}}{\varepsilon^{3\nu+2(p+1)}c_{p,d}^{2\nu+4(p+1)}(p+1)!^{4\nu}}\right]^\frac{1}{d\nu+2(p+1)(d+2\nu)}, \\
N&\asymp \left[\frac{B_\nu^{2d(p+1)}d^{2d\nu(p+1)}m^{d\nu+2(p+1)(d+\nu)}}{\varepsilon^{3d\nu+2(p+1)(2d+3\nu)}c_{p,d}^{d\nu+2d(p+1)}(p+1)!^{2d\nu}}\right]^\frac{1}{d\nu+2(p+1)(d+2\nu)}, \\
N_0&\asymp Nc_{p,d}\asymp \left[\frac{B_\nu^{2d(p+1)}c_{p,d}^{4\nu(p+1)}d^{2d\nu(p+1)}m^{d\nu+2(p+1)(d+\nu)}}{\varepsilon^{3d\nu+2(p+1)(2d+3\nu)}(p+1)!^{2d\nu}}\right]^\frac{1}{d\nu+2(p+1)(d+2\nu)}, \\
 R&\asymp \left[\frac{B_\nu^{d+4(p+1)}(p+1)!^{2d}m^{2(p+1)}}{\varepsilon^{2(p+1)-d}c_{p,d}^{4(p+1)}d^{2d(p+1)}}\right]^\frac{1}{d\nu+2(p+1)(d+2\nu)},
\end{align}
provided that $2(p+1)>d$ and $\nu>\frac{2(p+1)}{2(p+1)-d}$.
\footnote{Compare with footnote~\ref{ft:compl_strong} on page~\pageref{ft:compl_strong}.}
Thus, we have for the complexity
\begin{align}
\notag
\mathcal{C}_{ser}&\asymp JQ^dmc_{p,d}^2N\asymp JQ^dmc_{p,d}N_0\\
\label{compl_strong_ser}
&\asymp\left[\frac{B_\nu^{6d(p+1)}c_{p,d}^{2(p+1)(4\nu-d)-d\nu}d^{6d\nu(p+1)}m^{3d\nu+6(p+1)(d+\nu)}}{\varepsilon^{7d\nu+2(p+1)(4d+5\nu)}(p+1)!^{6d\nu}}\right]^\frac{1}{d\nu+2(p+1)(d+2\nu)}.
\end{align}
\end{theorem}

\begin{remark}
(i) Complexity estimate~\eqref{compl_strong_ser}
shows that one cannot go beyond the complexity order
$\varepsilon^{-2.5}$ in this case, no matter how large $p,\nu$ are. This is mainly due to the factor $J$ within the third constraint in~\eqref{mse_strong_constr_ser} which does not arise in~\eqref{mse_strong_constr}.

\smallskip
(ii) Similarly to Section~\ref{sec:compl_strong}, we would have obtained the same complexity even when we used a control variate 
with a higher variance order $\Delta^K$ for some $K>1$.

\smallskip
(iii) When comparing~\eqref{compl_strong_ser} with~\eqref{compl_strong}, one clearly sees that~\eqref{compl_strong} always achieves a better complexity for $\nu>\frac{2(p+1)}{2(p+1)-d}$ (in terms of $\varepsilon$).

\smallskip
(iv) Furthermore, also from the pure computational point
of view it is preferable to consider the integral approach rather than the series approach, even though the control variates $M_{\Delta,T}^{{ser},1}$ and $M_{\Delta,T}^{{int},1}$ are theoretically equivalent
(recall Theorem~\ref{th:13022017a1}).
This is mainly due to the factor $\Delta_j W^i$ in $a_{j,e_i}$ (see~\eqref{eq:wiener_pcv}), which is independent of $X_{\Delta,(j-1)\Delta}$ and has zero expectation and thus may lead to poor regression results (cf.~``RCV approach'' in~\cite{belomestny2016variance}).
Regarding the integral approach,
such a destabilising factor is not present in~$g_{j,k}$
(see~\eqref{def_g_strong}).
\end{remark}

\section{Numerical results}
\label{sec:5}

In this section, we consider the Euler scheme
and compare the numerical performance 
of the SMC, MLMC, series and integral approaches.
For simplicity we implemented a global regression
(i.e.\ the one without truncation and partitioning).
Regarding the choice of basis functions,
we use in both series and integral approaches
the same polynomials $\psi(x)=\prod_{k=1}^dx_k^{l_k}$,
where $l_1,\ldots l_d\in\left\{0,1,\ldots,p\right\}$
and $\sum_{k=1}^dl_k\leq p$.
In addition to the polynomials,
we consider the function $f$ as a basis function.
Hence, we have overall $\binom{p+d}{d}+1$
basis functions in each regression.
As for the MLMC approach,
we use the same simulation results
as in~\cite{belomestny2016variance}.

The following results are based on program codes
written and vectorised in \mbox{MATLAB}
and running on a Linux 64-bit operating system.

\subsection{One-dimensional example}
\label{num_1d}
\begin{figure}[htb!]
\includegraphics[width=0.49\textwidth]{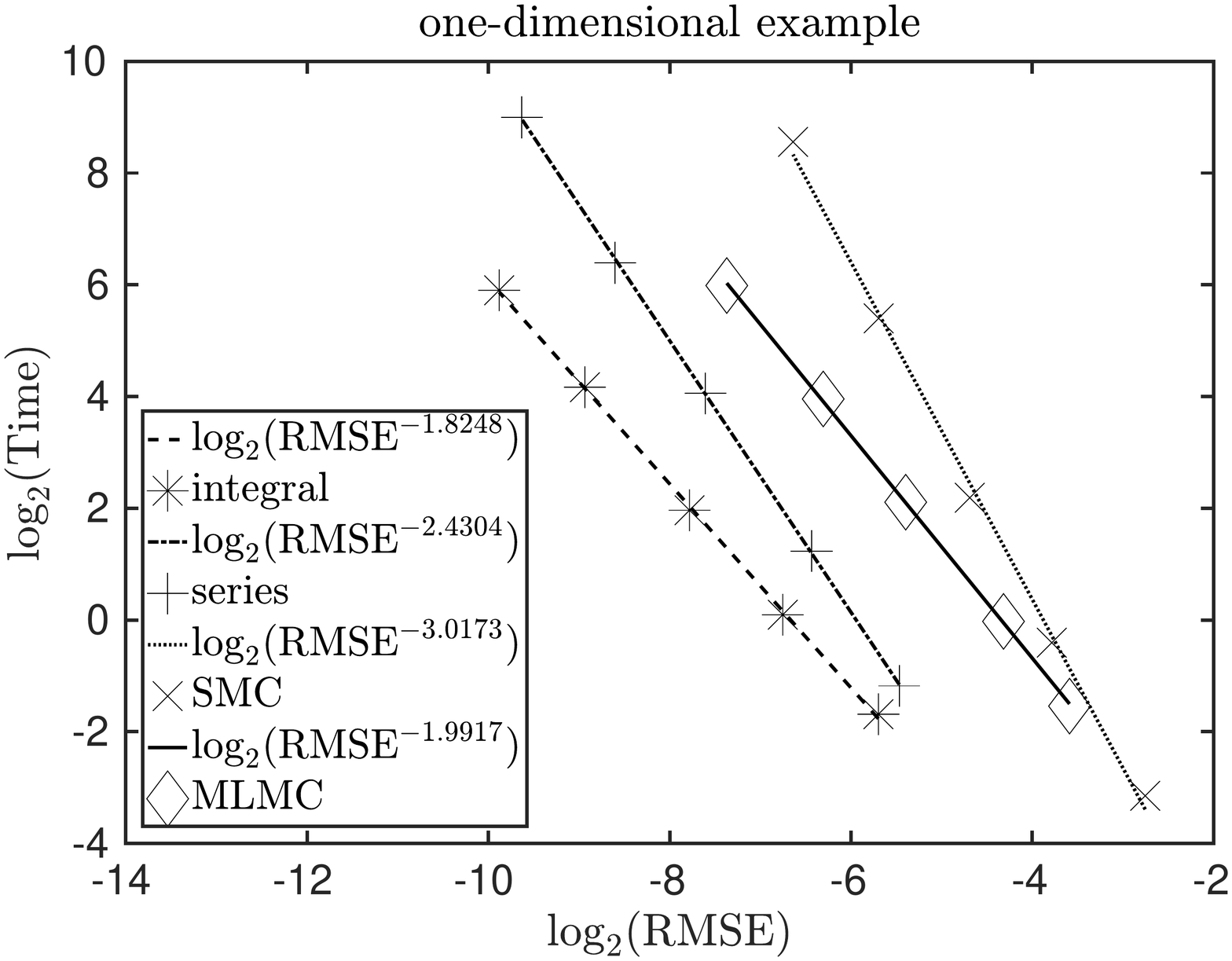}
\includegraphics[width=0.49\textwidth]{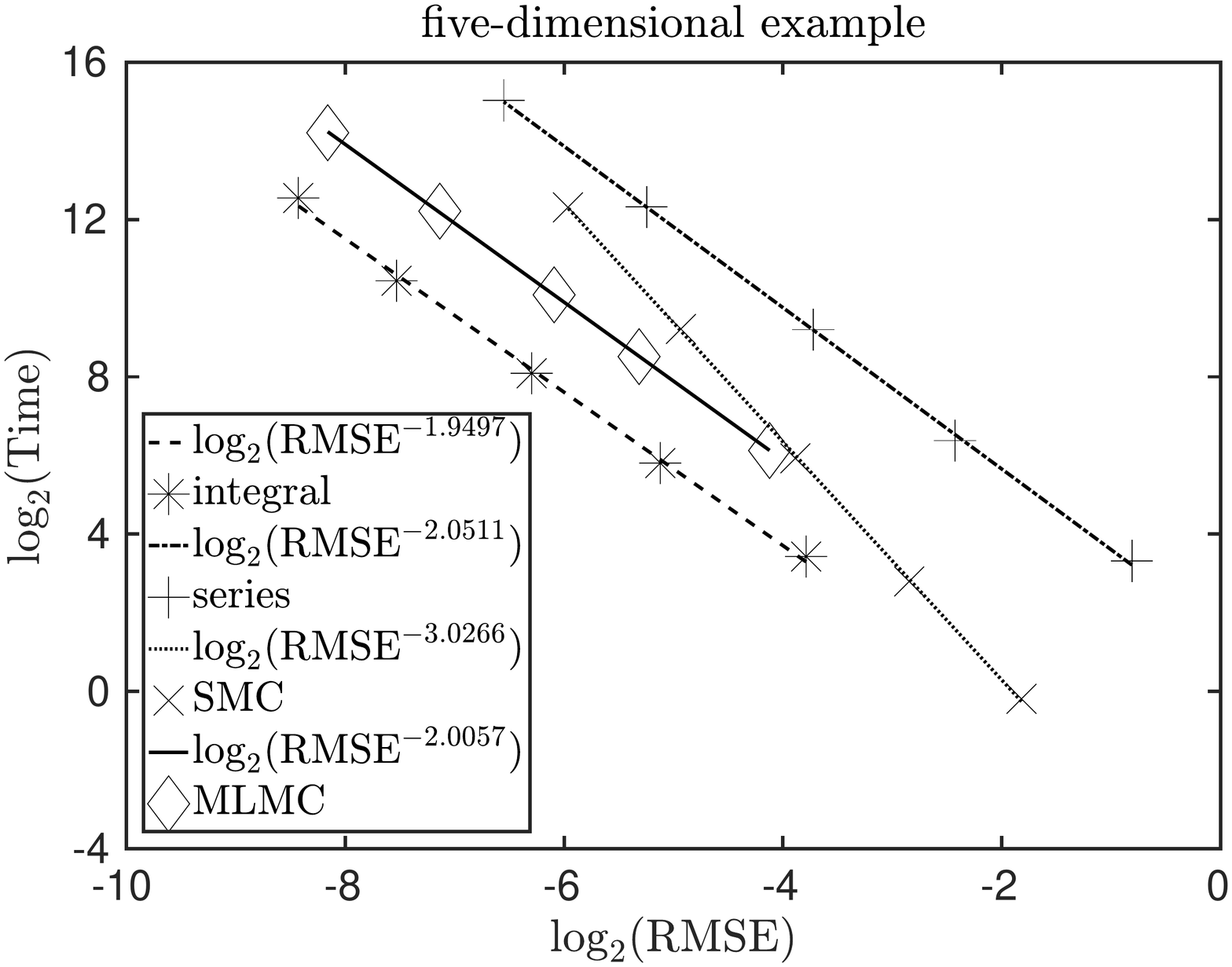}
\caption{Numerical complexities of the integral, series, SMC and MLMC approaches in the one- and five-dimensional cases.}
\label{compld}
\end{figure}

Here $d=m=1$.
We consider the following SDE (cf.~\cite{belomestny2016variance})
\begin{align}
\label{eq:sde1}
dX_t=&-\frac{1}{2}\tanh\left(X_t\right)\sech^2\left(X_t\right)dt+
\sech\left(X_t\right)dW_t,\quad X_0=0,
\end{align}
for $t\in\left[0,1\right]$, where $\sech(x)\doteq\frac{1}{\cosh(x)}$. This SDE has an exact solution \(X_t=\arsinh\left(W_t\right).\)
Furthermore, we consider the functional
$f(x)=\sech(x)+15\arctan(x)$, that is, we have
\begin{align}
\label{eq:0110a3}
\EE\left[f\left(X_1\right)\right]=
\EE\left[\sech\left(\arsinh\left(W_1\right)\right)\right]=
\EE\left[\frac{1}{\sqrt{1+W_1^2}}\right]\approx 0.789640.
\end{align}
We choose $p=3$ (that is, $5$~basis functions) and,
for each $\varepsilon=2^{-i}$, $i\in\left\{2,3,4,5,6\right\}$,
we set the parameters $J$, $N$ and \(N_0\) as follows
(compare with the formulas in Section~\ref{sec:4}
for $\nu\to\infty$, $\lim_{\nu\to\infty}B_\nu=1$ and ignore the log-terms for the integral approach):
\begin{align*}
J=\left\lceil \varepsilon^{-1}\right\rceil,\quad N=256\cdot\left\{\begin{array}{ll}\lceil 0.6342\cdot\varepsilon^{-1.0588}\rceil & \textnormal{integral approach}, \\ \lceil 0.6342\cdot\varepsilon^{-1.5882}\rceil & \textnormal{series approach},\end{array}\right. \\
N_0=256\cdot\left\{\begin{array}{ll}\lceil 2.5367\cdot\varepsilon^{-1.0588}\rceil & \textnormal{integral approach}, \\ \lceil 2.5367\cdot\varepsilon^{-1.5882}\rceil & \textnormal{series approach}.\end{array}\right.
\end{align*}
Regarding the SMC approach, the number of paths is set $N_0=256\cdot \varepsilon^{-2}$. The factor 256 is here for stability purposes.
As for the MLMC approach,
we set the initial number of paths in the first level ($l=0$)
equal to $10^3$ as well as the ``discretisation parameter''
$M=4$, which leads to time steps
of the length $\frac1{4^l}$ at level~$l$
(the notation here is as in~\cite{giles2008multilevel}).
Next we compute the numerical RMSE
(the exact value is known, see~\eqref{eq:0110a3}) by means of \(100\) independent repetitions of the algorithm.
As can be seen from left-hand side in Figure~\ref{compld},
the estimated numerical complexity is about $\text{RMSE}^{-1.82}$ for the integral approach, $\text{RMSE}^{-2.43}$ for the series approach, $\text{RMSE}^{-1.99}$ for the MLMC approach
and $\text{RMSE}^{-3.02}$ for the SMC approach,
which we get by regressing the log-time (logarithmic computing time of the whole algorithm in seconds) vs.\ log-RMSE.  
Thus, the 
complexity reduction works best
with the integral approach.

\subsection{Five-dimensional example}
Here $d=m=5$.
We consider the SDE (cf.~\cite{belomestny2016variance})
\begin{align}
\notag
dX_t^i&=-\sin\left(X_t^i\right)\cos^3\left(X_t^i\right)dt+\cos^2\left(X_t^i\right)dW_t^i,\quad X_0^i=0,\quad i\in\left\{1,2,3,4\right\},\\
\label{5d_sde}
dX_t^5&=\sum_{i=1}^4\left[-\frac{1}{2}\sin\left(X_t^i\right)\cos^2\left(X_t^i\right)dt+\cos\left(X_t^i\right)
dW_t^i\right]+dW_t^5,\quad X_0^5=0.
\end{align}
The solution of~\eqref{5d_sde} is given by
\begin{align*}
X_t^i&=\arctan\left(W_t^i\right),\quad i\in\left\{1,2,3,4\right\},\\
X_t^5&=\sum_{i=1}^4\operatorname{arsinh}\left(W_t^i\right)
+W_t^5.
\end{align*}
for $t\in\left[0,1\right]$. Further, we consider the functional 
\begin{align*}
f(x)=\cos\left(\sum_{i=1}^5x_i\right)-20\sum_{i=1}^4\sin\left(x_i\right),
\end{align*}
that is, we have
\begin{align*}
\mathbb{E}\left[f\left(X_1\right)\right]&=\left(\mathbb{E}\left[\cos\left(
\arctan\left(W_1^1\right)+\operatorname{arsinh}\left(W_1^1\right)\right)\right]\right)^4
\mathbb{E}\left[\cos\left(W_1^5\right)\right]
\approx 0.002069.
\end{align*}
We again choose $p=3$
(this now results in $57$ basis functions),
consider the same values of $\varepsilon$ as above (and, in addition, consider the values $\varepsilon=2^{-7}$ and $\varepsilon=2^{-8}$ for the SMC approach to obtain similar computing times as for the series and integral approaches). Moreover, we set (compare with the formulas in Section~\ref{sec:4}
for $\nu\to\infty$, $\lim_{\nu\to\infty}B_\nu=1$ and ignore the log-terms for the integral approach):
\begin{align*}
J=\left\lceil \varepsilon^{-1}\right\rceil,\quad N=\left\{\begin{array}{ll}\lceil 35.9733\cdot\varepsilon^{-1.2381}\rceil & \textnormal{integral approach}, \\ 4\cdot\lceil 4.9044\cdot\varepsilon^{-1.8571}\rceil & \textnormal{series approach},\end{array}\right. \\
N_0=\left\{\begin{array}{ll}\lceil 2014.5030\cdot\varepsilon^{-1.2381}\rceil & \textnormal{integral approach}, \\ 4\cdot\lceil 274.6480\cdot\varepsilon^{-1.8571}\rceil & \textnormal{series approach}.\end{array}\right.
\end{align*}
The number of paths for the SMC approach is again set $N_0=256\cdot \varepsilon^{-2}$. 
Regarding the MLMC approach, we again choose $M=4$,
but the initial number of paths in the first level is increased
to~$10^4$.
As in the one-dimensional case, we compute the numerical RMSE
by means of 100 independent repetitions of the algorithm.
Our empirical findings are illustrated on the right-hand side in Figure~\ref{compld}.
We observe the numerical complexity $\text{RMSE}^{-1.95}$ for the integral approach, $\text{RMSE}^{-2.05}$ for the series approach, $\text{RMSE}^{-2.01}$ for the MLMC approach
and $\text{RMSE}^{-3.03}$ for the SMC approach.
Even though here the complexity order of the series approach
is better than that of the SMC approach and close to that of MLMC approach,
the series approach is practically outperformed
by the other approaches
(see Figure~\ref{compld}; the multiplicative constant influencing
the computing time is obviously very big).
However, the integral approach remains numerically the best one
also in this five-dimensional example.

\section{Proofs}
\label{sec:proofs}

\subsection*{Proof of Theorem~\protect\ref{thm:ChaosDecompNum}}

Cf.~the proof of Theorem 2.1 in~\cite{belomestny2016variance}.

\subsection*{Proof of Theorem~\protect\ref{thm:IntegRepr}}
First of all, we derive 
\begin{align}
&\lim\limits_{t\nearrow t_j}u_\Delta(t,X_{\Delta,t_{j-1}},W_t-W_{t_{j-1}})\\
\notag
=&\lim\limits_{t\nearrow t_j}\EE\left[u_\Delta(t_j,\Phi_\Delta(x,y+z_j\sqrt{t_j-t}),0)\right]\left|\,_{x=X_{\Delta,(j-1)\Delta},\,y=W_t-W_{t_{j-1}}}\right.\\
\notag
=&u_\Delta(t_j,\Phi_\Delta(X_{\Delta,(j-1)\Delta},\Delta_jW),0)=u_\Delta(t_j,X_{\Delta,t_j},0).
\end{align}
By means of It\^o's lemma and the fact that $u_\Delta$ satisfies the heat equation
\begin{align}
\label{heat_u}
\frac{\partial u_\Delta}{\partial t}
+\frac12\sum_{i=1}^m\frac{\partial^2 u_\Delta}{\partial y_i^2}
=0
\end{align}
due to its relation to the normal distribution, we then obtain
\begin{align}
\label{repr_integ_proof}
&f(X_{\Delta,T})-\EE[f(X_{\Delta,T})]\\
\notag
=&u_\Delta(T,X_{\Delta,T},0)-u_\Delta(0,x_0,0)\\
\notag
=&\sum_{j=1}^J\left(u_\Delta(t_j,X_{\Delta,t_j},0)-u_\Delta(t_{j-1},X_{\Delta,t_{j-1}},0)\right)\\
\notag
=&\sum_{j=1}^J\lim\limits_{t\nearrow t_j}\left( u_\Delta(t,X_{\Delta,t_{j-1}},W_t-W_{t_{j-1}})-u_\Delta(t_{j-1},X_{\Delta,t_{j-1}},0)\right)\\
\notag
=&\sum_{j=1}^J\sum_{i=1}^m\lim\limits_{t\nearrow t_j}\,\int\limits_{t_{j-1}}^t\frac{\partial u_\Delta}{\partial y_i}(s,X_{\Delta,t_{j-1}},W_s-W_{t_{j-1}})\,dW_s^i\\
\notag
=&\sum_{j=1}^J\,\int\limits_{t_{j-1}}^{t_j}\nabla_y u_\Delta(s,X_{\Delta,t_{j-1}},W_s-W_{t_{j-1}})\,dW_s.
\end{align}
Next, let us derive a relation between $\nabla_y u_\Delta$ and $\nabla_x u_\Delta$. We have for $t\in[t_{j-1},t_j)$
\begin{align*}
\nabla_y u_\Delta(t,x,y)&=\nabla_y \EE[u_\Delta(t_j,\Phi_\Delta(x,y+z_j\sqrt{t_j-t}),0)]\\
&=\nabla_x\EE[u_\Delta(t_j,\Phi_\Delta(x,y+z_j\sqrt{t_j-t}),0)]\sigma(x).
\end{align*}
Thus, the term $\nabla_y u_\Delta(s,X_{\Delta,t_{j-1}},W_s-W_{t_{j-1}})$ in~\eqref{repr_integ_proof} takes the form
\begin{align}
\label{deriv_x_y}
\nabla_y u_\Delta(s,X_{\Delta,t_{j-1}},W_s-W_{t_{j-1}})=\EE[\nabla_xu_\Delta(t_j,X_{\Delta,t_j},0)\left|\,\mathcal{F}_s\right.]\sigma(X_{\Delta,t_{j-1}}).
\end{align}
Note that it holds
\begin{align*}
u_\Delta(t_j,x,0)=\EE[f(X_{t_j,x}(\Delta,T))],
\end{align*}
where we recall that $X_{t_j,x}(\Delta,t_l)$, for $l\ge j$, denotes the Euler discretisation starting at time $t_j$ in $x$ (analogous to $X_{s,x}(t)$ for the exact solution). Hence, we have for $\nabla_x u_\Delta$
\begin{align*}
\nabla_xu_\Delta(t_j,x,0)=\EE[\nabla f(X_{t_j,x}(\Delta,T))\delta X_{t_j,x}(\Delta,T)]
\end{align*}
or, in another form,
\begin{align*}
\nabla_xu_\Delta(t_j,X_{\Delta,t_j},0)=\EE\left[\nabla f(X_{t_j,X_{\Delta,t_j}}(\Delta,T))\delta X_{t_j,X_{\Delta,t_j}}(\Delta,T)\left|\,\mathcal{F}_{t_j}\right.\right],
\end{align*}
where  $\delta^i X^k_{t_j,x}(\Delta,t_l)\doteq\frac{\partial X^k_{t_j,x}(\Delta,t_l)}{\partial x_i}$ with $l\ge j$ and $i,k\in\left\{1,\ldots,d\right\}$.
We also notice at this point that
$X_{\Delta,t_l}=X_{0,x_0}(\Delta,t_l)$ and $\delta X_{\Delta,t_l}=\delta X_{0,x_0}(\Delta,t_l)$.

Let us define $\sigma_k(x)\doteq\left(\sigma_{k,1}(x),\ldots,\sigma_{k,m}(x)\right)^\top$ for $k\in\left\{1,\ldots,d\right\}$. Further, we denote with $\mathcal{J}_\mu\in\bbR^{d\times d}$, $\mathcal{J}_{\sigma_k}\in\bbR^{m\times d}$ the Jacobi matrices of the functions $\mu$, $\sigma_k$. Regarding the discretisation $\delta X_{\Delta,j\Delta}$ of $\delta X_t$ we can use, alternatively to~\eqref{deriv_discr_euler}, the matrix form
\begin{align}
\delta X_{\Delta, j\Delta} 
    = A_j\delta X_{\Delta,(j-1)\Delta}=A_jA_{j-1}\cdots A_1,
\end{align}
where 
\begin{align*}
A_k\doteq I_d+\mathcal{J}_\mu(X_{\Delta,(k-1)\Delta})\Delta+\begin{pmatrix}\Delta_k W^\top \mathcal{J}_{\sigma_1}(X_{\Delta,(k-1)\Delta})\\ \vdots \\ \Delta_k W^\top \mathcal{J}_{\sigma_d}(X_{\Delta,(k-1)\Delta})\end{pmatrix}.
\end{align*}
This gives us 
\begin{align*}
X_{t_j,X_{\Delta,t_j}}(\Delta,t_l)&=\Phi_\Delta(\cdots(\Phi_\Delta(X_{\Delta,t_j},\Delta_{j+1}W),\cdots,\Delta_{l}W)\\
&=\Phi_\Delta(\cdots(\Phi_\Delta(X_{\Delta,0},\Delta_{1}W),\cdots,\Delta_{l}W)=X_{\Delta,t_l},\\
\delta X_{t_j,X_{\Delta,t_j}}(\Delta,t_l)&=A_lA_{l-1}\cdots A_{j+1}=A_lA_{l-1}\cdots A_1\left(A_jA_{j-1}\cdots A_1\right)^{-1}\\
&=\delta X_{\Delta,t_l}\delta X_{\Delta,t_j}^{-1},
\end{align*}
where $\Phi_\Delta$ is defined through~\eqref{phi_euler}.
Finally, we obtain for $s\in\left[t_{j-1},t_j\right)$
\begin{align*}
\nabla_y u_\Delta(s,X_{\Delta,t_{j-1}},W_s-W_{t_{j-1}})&=\EE\left[\EE\left[\nabla f(X_{\Delta,T})\delta X_{\Delta,T}\delta X_{\Delta,t_j}^{-1}\left|\,\mathcal{F}_{t_j}\right.\right]\left|\,\mathcal{F}_s\right.\right]\sigma(X_{\Delta,t_{j-1}})
\\
&=\EE\left[\nabla f(X_{\Delta,T})\delta X_{\Delta,T}\delta X_{\Delta,t_j}^{-1}\left|\,\mathcal{F}_s\right. \right]\sigma(X_{\Delta,t_{j-1}}).
\end{align*}

\subsection*{Proof of Theorem~\protect\ref{relation_eul}}
Below we simply write $u_{\Delta,t_{j-1}}$ rather than $u_{\Delta}(t_{j-1},X_{\Delta,t_{j-1}},0)$. Let us consider the Taylor expansion for $\frac{\partial}{\partial y_r}u_\Delta(t,X_{\Delta,t_{j-1}},W_t-W_{t_{j-1}})$ of order $K\in\NN_0$ around $(t_{j-1},X_{\Delta,t_{j-1}},0$), with $r\in\left\{1,\ldots,m\right\}$, that is,
for $t\in[t_{j-1},t_j)$, we set
\begin{align}
\label{taylor}
T_{j,r}^K(t)\doteq \sum_{\left|\alpha\right|\le K}\frac{D^\alpha\left(\frac{\partial}{\partial y_r}u_{\Delta,t_{j-1}}\right)}{\alpha_1!\cdots\alpha_{m+1}!}(t-t_{j-1})^{\alpha_1}(W_t^1-W_{t_{j-1}}^1)^{\alpha_2}\cdots(W_t^m-W_{t_{j-1}}^m)^{\alpha_{m+1}},
\end{align}
where $\alpha\in\NN_0^{m+1}$ and $D^\alpha \left(\frac{\partial}{\partial y_m}u_{\Delta,t_{j-1}}\right)= \frac{\partial^{\left|\alpha\right|}\left(\frac{\partial}{\partial y_m}u_{\Delta,t_{j-1}}\right)}{\partial t^{\alpha_1}\partial y_1^{\alpha_2}\cdots\partial y_m^{\alpha_{m+1}}}$.
Via Taylor's theorem we obtain
\begin{align*}
&\frac{\partial}{\partial y_r}u_\Delta(t,X_{\Delta,t_{j-1}},W_t-W_{t_{j-1}})-T_{j,r}^K(t)\\
=&\sum_{\left|\alpha\right|=K+1}\left[\frac{(K+1)!}{\alpha_1!\cdots\alpha_{m+1}!}\int\limits_0^1(1-z)^KD^\alpha\left(\frac{\partial}{\partial y_r}u_\Delta(t_{j-1}+z(t-t_{j-1}),X_{\Delta,t_{j-1}},z(W_t-W_{t_{j-1}}))\right)\,dz\right.\\
&\phantom{\sum_{\left|\alpha\right|=K+1}\left[\right.}\left.\cdot (t-t_{j-1})^{\alpha_1}(W_t^1-W_{t_{j-1}}^1)^{\alpha_2}\cdots(W_t^m-W_{t_{j-1}}^m)^{\alpha_{m+1}}\right].
\end{align*}
Provided that~\eqref{deriv_bound} holds, we get
\begin{align*}
&\Var\left[\sum_{j=1}^J\sum_{r=1}^m\,\int\limits_{t_{j-1}}^{t_j}\left(\frac{\partial}{\partial y_r} u_\Delta(t,X_{\Delta,t_{j-1}},W_t-W_{t_{j-1}})-T_{j,r}^K(t)\right)\,dW_t^r\right]\\
&=\sum_{j=1}^J\sum_{r=1}^m\,\int\limits_{t_{j-1}}^{t_j}\EE\left[\left(\frac{\partial}{\partial y_r} u_\Delta(t,X_{\Delta,t_{j-1}},W_t-W_{t_{j-1}})-T_{j,r}^K(t)\right)^2\right]\,dt\\
\lesssim &C^{2(K+1)}\sum_{j=1}^J\sum_{\left|\alpha\right|=K+1}\,\int\limits_{t_{j-1}}^{t_j}\EE\left[(t-t_{j-1})^{2\alpha_1}(W_t^1-W_{t_{j-1}}^1)^{2\alpha_2}\cdots(W_t^m-W_{t_{j-1}}^m)^{2\alpha_{m+1}}\right]\,dt\\
\lesssim & (C^2\Delta)^{K+1}\stackrel{K\to\infty}{\longrightarrow}0,
\end{align*}
and thus $T_{j,r}^K$ converges for $K\to\infty$
in $L^2(\Omega\times[0,T])$
to $\frac{\partial u_\Delta}{\partial y_r}(t,X_{\Delta,t_{j-1}},W_t-W_{t_{j-1}})$.
Moreover, due to~\eqref{heat_u}, the limit of $T_{j,r}^K$ simplifies to (cf.~\eqref{taylor})
\begin{align*}
&\frac{\partial u_{\Delta,t_{j-1}}}{\partial y_r}+\sum_{i=1}^m\frac{\partial^2 u_{\Delta,t_{j-1}}}{\partial y_r\partial y_i}(W_t^i-W_{t_{j-1}}^i)\\
&+\frac{1}{2}\sum_{i=1}^m\frac{\partial^3 u_{\Delta,t_{j-1}}}{\partial y_r\partial y_i^2}((W_t^i-W_{t_{j-1}}^i)^2-(t-t_{j-1}))+\sum_{\substack{i_1,i_2=1\\i_1<i_2}}^m\frac{\partial^3 u_{\Delta,t_{j-1}}}{\partial y_r\partial y_{i_1}\partial y_{i_2}}(W_t^{i_1}-W_{t_{j-1}}^{i_1})(W_t^{i_2}-W_{t_{j-1}}^{i_2}) \\
&+\left[\frac{1}{6}\sum_{i=1}^m\frac{\partial^4 u_{\Delta,t_{j-1}}}{\partial y_r\partial y_i^3}((W_t^i-W_{t_{j-1}}^i)^3-3(W_t^i-W_{t_{j-1}}^i)(t-t_{j-1}))\right.\\
\notag
&\phantom{+\left[\right.}+\frac{1}{2}\sum_{\substack{i_1,i_2=1\\i_1<i_2}}^m\frac{\partial^4 u_{\Delta,t_{j-1}}}{\partial y_r\partial y_{i_1}^2\partial y_{i_2}}((W_t^{i_1}-W_{t_{j-1}}^{i_1})^2-(t-t_{j-1}))(W_t^{i_2}-W_{t_{j-1}}^{i_2}) \\
\notag
&\left.\phantom{+\left[\right.}+\sum_{\substack{i_1,i_2,i_3=1\\i_1<i_2<i_3}}^m\frac{\partial^4 u_{\Delta,t_{j-1}}}{\partial y_r\partial y_{i_1}\partial y_{i_2}\partial y_{i_3}}(W_t^{i_1}-W_{t_{j-1}}^{i_1})(W_t^{i_2}-W_{t_{j-1}}^{i_2})(W_t^{i_3}-W_{t_{j-1}}^{i_3})\right] \\
\notag
&+...\\
=&\sum_{l=1}^{\infty}(t-t_{j-1})^\frac{l-1}{2}\sum_{\substack{k\in\NN_0^m
\\ \sum_{i=1}^mk_i=l-1}}\frac{\partial ^l u_{\Delta,t_{j-1}}}{\partial y_r \partial y_1^{k_1}\cdots \partial y_m^{k_m}}\prod_{i=1}^m\frac{H_{k_i}\left(\frac{W_t^i-W_{t_{j-1}}^i}{\sqrt{t-t_{j-1}}}\right)}{\sqrt{k_i!}}.
\end{align*}
To compute the stochastic integral
\begin{align*}
&\int\limits_{t_{j-1}}^{t_j}\nabla_y u_\Delta(t,X_{\Delta,t_{j-1}},W_t-W_{t_{j-1}})\,dW_t\\
=&\sum_{l=1}^{\infty}\sum_{r=1}^m\int\limits_{t_{j-1}}^{t_j}(t-t_{j-1})^\frac{l-1}{2}\sum_{\substack{k\in\NN_0^m
\\ \sum_{i=1}^mk_i=l-1}}\frac{\partial ^l u_{\Delta,t_{j-1}}}{\partial y_r \partial y_1^{k_1}\cdots \partial y_m^{k_m}}\prod_{i=1}^m\frac{H_{k_i}\left(\frac{W_t^i-W_{t_{j-1}}^i}{\sqrt{t-t_{j-1}}}\right)}{\sqrt{k_i!}}\,dW_t^r,
\end{align*}
we apply It\^o's lemma w.r.t.\ the functions $F_k(t,y_1,\ldots,y_m)\doteq t^{l/2}\prod_{i=1}^m\frac{H_{k_i}\left(\frac{y_i}{\sqrt{t}}\right)}{\sqrt{k_i!}}$, where $\sum_{i=1}^mk_i=l$.
Thus, we obtain
\begin{align}
\label{ito_strong}
&dF_k(t-t_{j-1},W_t^1-W_{t_{j-1}}^1,\ldots,W_t^m-W_{t_{j-1}}^m)\\
\notag
=&(t-t_{j-1})^\frac{l-1}{2}\sum_{r=1}^m\frac{H_{k_r-1}\left(\frac{W_t^r-W_{t_{j-1}}^r}{\sqrt{t-t_{j-1}}}\right)}{\sqrt{(k_r-1)!}}\prod_{\substack{i=1\\ i\neq r}}^m\frac{H_{k_i}\left(\frac{W_t^i-W_{t_{j-1}}^i}{\sqrt{t-t_{j-1}}}\right)}{\sqrt{k_i!}}\,dW_t^r.
\end{align}
This gives us finally
\begin{align*}
&\int\limits_{t_{j-1}}^{t_j}\nabla_y u_\Delta(t,X_{\Delta,t_{j-1}},W_t-W_{t_{j-1}})\,dW_t\\
=&\sum_{l=1}^\infty\Delta^{l/2}\sum_{\substack{k\in\NN_0^m\\ \sum_{i=1}^mk_i=l}}\frac{\partial^lu_{\Delta}(t_{j-1},X_{\Delta,t_{j-1}},0)}{\partial y_1^{k_1}\cdots\partial y_m^{k_m}}\prod_{i=1}^m\frac{H_{k_i}\left(\frac{\Delta_jW^i}{\sqrt{\Delta}}\right)}{\sqrt{k_i!}}.
\end{align*}

\subsection*{Proof of Theorem~\protect\ref{th:13022017a1}}
We define the (random) function $G_{l,j}(x)$ for $J\ge l\ge j\ge 0$, $x\in\RR^d$, as follows
\begin{align}
\label{strong:def_G}
G_{l,j}(x)&=\Phi_{\Delta,l}\circ\Phi_{\Delta,l-1}\circ\ldots\circ
\Phi_{\Delta,{j+1}}(x),\quad l>j,\\
\notag
G_{l,j}(x)&=x,\quad l=j,
\end{align}
where $\Phi_{\Delta,l}(x)\doteq\Phi_{\Delta}\left(x,\Delta_l W\right)$ for $l=1,\ldots,J$. 
Note that  it holds
\begin{align}
\label{strong:rel_uG2}
u_\Delta(t_j,x,0)=\EE\left[f(G_{J,j}(x))\right].
\end{align}
Similar to $G$ we define the function $\tilde G_{j}(x,z)$, $0\le j< J$, $x\in\RR^d$, $z:=(z_1,\ldots,z_{J-j})\in\RR^{m\times (J-j)}$, $z_l:=(z_l^1,\ldots,z_l^m)^\top\in\RR^m$ for $l=1,\ldots,J-j$, as follows
$$
\tilde G_{j}(x,z):=\tilde\Phi_{\Delta,z_{J-j}}\circ\ldots\circ\tilde\Phi_{\Delta,z_{1}}(x),
$$
where $\tilde\Phi_{\Delta,z_l}(x)\doteq\Phi_{\Delta}\left(x,z_l\sqrt\Delta\right)$. Note that $G$ and $\tilde G$ are in the following relation
\begin{align}
\label{strong:rel_G}
G_{J,j}(x)=\tilde G_j\left(x,\frac{1}{\sqrt{\Delta}}\begin{pmatrix}\Delta_{j+1} W,\Delta_{j+2}W,\ldots,\Delta_J W\end{pmatrix}\right),\quad j<J.
\end{align}

Let us represent $\sqrt{\Delta}\frac{\partial}{\partial y_i}u_\Delta(t_{j-1},x,0)$, where $j\in\left\{1,\ldots,J\right\}$ and $i\in\left\{1,\ldots,m\right\}$, as a $(J-j+1)m$-dimensional integral, that is (cf.~\eqref{strong:rel_G})
\begin{align*}
&\sqrt{\Delta}\frac{\partial}{\partial y_i}u_\Delta(t_{j-1},x,0)=\sqrt{\Delta}\frac{\partial}{\partial y_i}\EE\left[f\left(G_{J,j}\left(\Phi_\Delta\left(x,\Delta_j W+y\right)\right)\right)\right]\left.\right|_{y=0_m}\\
=&\int\limits_{\RR^{(J-j+1)m}}\sqrt{\Delta}\frac{\partial}{\partial y_i}
\left[
f\left(\tilde G_{j-1}\left(x,\left(z_1+\frac{y}{\sqrt{\Delta}},z_2,\ldots,z_{J-j+1}\right)\right) \right)
\right]
\varphi_{(J-j+1)m}(z)\, dz\left|_{y=0_m}\right.,
\end{align*}
where
$\varphi_{(J-j+1)m}$ denotes the $(J-j+1)m$-dimensional standard normal density function. Since it holds
\begin{align*}
&\sqrt{\Delta}\frac{\partial}{\partial y_i}
\left[
f\left(\tilde G_{j-1}\left(x,\left(z_1+\frac{y}{\sqrt{\Delta}},z_2,\ldots,z_{J-j+1}\right)\right) \right)
\right]
\\
=&\frac{\partial}{\partial z_1^i}
\left[
f\left(\tilde G_{j-1}\left(x,\left(z_1+\frac{y}{\sqrt{\Delta}},z_2,\ldots,z_{J-j+1}\right)\right) \right)
\right],
\end{align*}
we obtain via integration by parts
\begin{align*}
&\sqrt{\Delta}\frac{\partial}{\partial y_i}u_\Delta(t_{j-1},x,0)\\
=&\int\limits_{\RR^{(J-j+1)m}}\frac{\partial}{\partial z_1^i}
\left[
f\left(\tilde G_{j-1}\left(x,z\right) \right)
\right]
\varphi_{(J-j+1)m}(z)\, dz\\
=&-\int\limits_{\RR^{(J-j+1)m}}f\left(\tilde G_{j-1}\left(x,z\right) \right)\frac{\partial}{\partial z_1^i}\varphi_{(J-j+1)m}(z)\, dz\\
=&\int\limits_{\RR^{(J-j+1)m}}f\left(\tilde G_{j-1}\left(x,z\right) \right)z_1^i\varphi_{(J-j+1)m}(z)\, dz\\
=&\EE\left[f(G_{J,j-1}(x))\frac{\Delta_j W^i}{\sqrt{\Delta}}\right]=\EE\left[f(X_{\Delta,T})\frac{\Delta_jW^i}{\sqrt{\Delta}}\left|X_{\Delta,(j-1)\Delta}=x\right.\right]=a_{j,e_i}(x).
\end{align*}
We finally remark that we have
only the integral term in the integration by parts above
because the function
$z_1\mapsto f(\tilde G_{j-1}(x,z))\varphi_{(J-j+1)m}(z)$
is integrable over $\bbR$
w.r.t.\ the Lebesgue measure.

\subsection*{Proof of Theorem~\protect\ref{var_ord_eul}}

Via Taylor's theorem we get
\begin{align}
\notag
&\frac{\partial u_{\Delta}(t,X_{\Delta,t_{j-1}},W_t-W_{t_{j-1}})}{\partial y_i}\\
\notag
=&\frac{\partial u_{\Delta}(t_{j-1},X_{\Delta,t_{j-1}},0)}{\partial y_i}+(t-t_{j-1})\int\limits_0^1\frac{\partial^2 u_{\Delta}(t_{j-1}+z(t-t_{j-1}),X_{\Delta,t_{j-1}},z(W_t-W_{t_{j-1}}))}{\partial y_i\partial t}\,dz\\
\label{taylor_u}
&+\sum_{r=1}^m(W_t^r-W_{t_{j-1}}^r)\int\limits_0^1\frac{\partial^2 u_{\Delta}(t_{j-1}+z(t-t_{j-1}),X_{\Delta,t_{j-1}},z(W_t-W_{t_{j-1}}))}{\partial y_i\partial y_r}\,dz.
\end{align}
Due to~\eqref{heat_u}, \eqref{taylor_u}~simplifies to
\begin{align*}
&\frac{\partial u_{\Delta}(t,X_{\Delta,t_{j-1}},W_t-W_{t_{j-1}})}{\partial y_i}\\
=&\frac{\partial u_{\Delta}(t_{j-1},X_{\Delta,t_{j-1}},0)}{\partial y_i}-\frac12(t-t_{j-1})\int\limits_0^1\sum_{r=1}^m\frac{\partial^3 u_{\Delta}(t_{j-1}+z(t-t_{j-1}),X_{\Delta,t_{j-1}},z(W_t-W_{t_{j-1}}))}{\partial y_i\partial y_r^2}\,dz\\
&+\sum_{r=1}^m(W_t^r-W_{t_{j-1}}^r)\int\limits_0^1\frac{\partial^2 u_{\Delta}(t_{j-1}+z(t-t_{j-1}),X_{\Delta,t_{j-1}},z(W_t-W_{t_{j-1}}))}{\partial y_i\partial y_r}\,dz.
\end{align*}
Provided that the second and third derivatives of $u_\Delta$ w.r.t. $y$ are bounded, we have
\begin{align*}
&\Var\left[\int\limits_{t_{j-1}}^{t_j}\frac{\partial u_{\Delta}(t,X_{\Delta,t_{j-1}},W_t-W_{t_{j-1}})}{\partial y_i}\,dW_{t}^i-\frac{\partial u_{\Delta}(t_{j-1},X_{\Delta,t_{j-1}},0)}{\partial y_i}\Delta_jW^i\right]\\
=&\int\limits_{t_{j-1}}^{t_j}\EE\left[\left(\int\limits_0^1\sum_{r=1}^m\left((W_t^r-W_{t_{j-1}}^r)\frac{\partial^2 u_{\Delta}(t_{j-1}+z(t-t_{j-1}),X_{\Delta,t_{j-1}},z(W_t-W_{t_{j-1}}))}{\partial y_i\partial y_r}\right.\right.\right.\\
&\phantom{\int\limits_{t_{j-1}}^{t_j}\EE\left[\left(\int\limits_0^1\sum_{r=1}^m\left(\right.\right.\right.}\left.\left.\left.-\frac12(t-t_{j-1})\frac{\partial^3 u_{\Delta}(t_{j-1}+z(t-t_{j-1}),X_{\Delta,t_{j-1}},z(W_t-W_{t_{j-1}}))}{\partial y_i\partial y_r^2} \right)\,dz\right)^2\right]\,dt\\
\lesssim&
\sum_{r=1}^m\,\int\limits_{t_{j-1}}^{t_j}\EE\left[(W_t^r-W_{t_{j-1}}^r)^2+(t-t_{j-1})^2\right]\,dt
\lesssim\Delta^2.
\end{align*}
Thus, we finally obtain
\begin{align*}
\Var\left[f(X_{\Delta,T})-M_{\Delta,T}^{int,1}\right]\lesssim 
\Delta.
\end{align*}

\subsection*{Proof of Theorem~\protect\ref{func_assump}}
We start the calculations,
which will lead to the proof of part~(ii).
At some point we will get the proof of part~(i)
as a by-product.

In this proof we will use the shorthand notation
$\xi_k:=\Delta_k W$, $k\in\{1,\ldots,J\}$.
For $j\in\{0,\ldots,J-1\}$, we have
\[
u_{\Delta}(t_{j},x,y)=\EE\left[f\left(\Phi_{\Delta,J}\circ\Phi_{\Delta,J-1}\circ\ldots
\circ\Phi_{\Delta,j+2}
\circ\Phi_{\Delta}(x,y+\xi_{j+1})\right)\right],
\]
where $\Phi_{\Delta,k}(x)\doteq\Phi_{\Delta}(x,\xi_k)$.
Denote, for $k>j$,
\[
G_{k,j}(x,y)\doteq\Phi_{\Delta,k}\circ\Phi_{\Delta,k-1}\circ\ldots
\circ\Phi_{\Delta,j+2}
\circ\Phi_{\Delta}(x,y+\xi_{j+1}).
\]
Assume that for any $n\in\NN$, $l\in\left\{1,\ldots,d\right\}$, $\alpha\in\NN_0^d$,
\begin{align}
\label{Phi_assump_1}
\left|\EE\left[\left.\left(D^\alpha\Phi{}_{\Delta,k+1}^{l}(G_{k,j}(x,y))\right)^{n}\right|\mathcal{F}_{k}\right]\right|\leq\begin{cases}
(1+A_{n,l}\Delta), & \beta=\alpha_l=1\\
B_{n,l,\alpha}\Delta, & (\beta>1)\vee (\alpha_l\ne 1)
\end{cases}
\end{align}
with probability one for $\beta=\left|\alpha\right|\in\NN$ and some constants $A_{n,l}>0$, $B_{n,l,\alpha}>0$.
We recall the notation $D^\alpha f(x)= \frac{\partial^{\left|\alpha\right|}f(x)}{\partial x_1^{\alpha_1}\cdots\partial x_d^{\alpha_d}}$,
which was used here.
Clearly, for the Euler scheme~\eqref{phi_euler},
condition~\eqref{Phi_assump_1} is satisfied if all the derivatives of order $\beta$ for $\mu_k,\sigma_{ki}$, $k\in\left\{1,\ldots,d\right\}$, $i\in\left\{1,\ldots,m\right\}$, are bounded.
Moreover, suppose that for any $n_1,n_2\in\NN$, $l\in\left\{1,\ldots,d\right\}$, $\alpha_1,\alpha_2\in\NN_0^d$, with $\beta_1=\left|\alpha_1\right|>0$, $\beta_2=\left|\alpha_2\right|>0$, $(\beta_1>1)\vee (\beta_2>1)\vee ((\alpha_1)_{l}\neq 1)\vee((\alpha_2)_{l}\neq 1)$,
\begin{align}
\label{Phi_assump_2}
\left|\EE\left[\left.\left(D^{\alpha_1}\Phi{}_{\Delta,k+1}^{l}(G_{k,j}(x,y))\right)^{n_1}\left(D^{\alpha_2}\Phi{}_{\Delta,k+1}^{l}(G_{k,j}(x,y))\right)^{n_2}\right|\mathcal{F}_{k}\right]\right|\leq C_{n_1,n_2,l,\alpha_1,\alpha_2}\Delta
\end{align}
for some constant $C_{n_1,n_2,l,\alpha_1,\alpha_2}>0$.
Again, for the Euler scheme~\eqref{phi_euler},
condition~\eqref{Phi_assump_2} is satisfied if all the derivatives of orders $\beta_1$ and $\beta_2$ for $\mu_k,\sigma_{ki}$ are bounded.

We have for some $i\in\left\{1,\ldots,m\right\}$ and $l\in\left\{1,\ldots,d\right\}$
\[
\frac{\partial}{\partial y_i}G_{k+1,j}^l(x,y)=\sum_{s=1}^d\frac{\partial}{\partial x_s}\Phi_{\Delta,k+1}^{l}(G_{k,j}(x,y))\frac{\partial}{\partial y_i}G_{k,j}^s(x,y)
\]
and $\frac{\partial}{\partial y_i}G_{j+1,j}^s(x,y)=\frac{\partial}{\partial y_i}\Phi_{\Delta}^s(x,y+\xi_{j+1})$.
Hence 
\begin{eqnarray*}
&&\EE\left[\left(\frac{\partial}{\partial y_i}G_{k+1,j}^l(x,y)\right)^{2}\right]\\
&\leq&\EE\left[\left(\frac{\partial}{\partial x_l}\Phi_{\Delta,k+1}^{l}(G_{k,j}(x,y))\frac{\partial}{\partial y_i}G_{k,j}^l(x,y)\right)^2\right.\\
&&\phantom{\EE\left[\right.}+\sum_{s\neq l}\left\{2\frac{\partial}{\partial x_l}\Phi_{\Delta,k+1}^{l}(G_{k,j}(x,y))\frac{\partial}{\partial x_s}\Phi_{\Delta,k+1}^{l}(G_{k,j}(x,y))\frac{\partial}{\partial y_i}G_{k,j}^l(x,y)\frac{\partial}{\partial y_i}G_{k,j}^s(x,y)\right.\\
&&\phantom{\EE\left[+\sum_{s\neq l}\left\{\right.\right.}\left.\left.+(d-1)\left(\frac{\partial}{\partial x_s}\Phi_{\Delta,k+1}^{l}(G_{k,j}(x,y))\frac{\partial}{\partial y_i}G_{k,j}^s(x,y)\right)^2\right\}\right]
\end{eqnarray*}
Denote 
\begin{align}
\label{eq:14022017a1}
\rho_{k+1,n}^{i,s}=\EE\left[\left(\frac{\partial}{\partial y_i}G_{k+1,j}^s(x,y)\right)^{n}\right],
\end{align}
then, due to $2ab\le a^2+b^2$, we get for $k=j+1,\ldots,J-1$,
\begin{eqnarray*}
\rho_{k+1,2}^{i,l} & \leq & (1+A_{2,l}\Delta)\rho_{k,2}^{i,l}
+\sum_{s\neq l}\left\{C_{1,1,l,e_l,e_s}\Delta(\rho_{k,2}^{i,l}+\rho_{k,2}^{i,s})+(d-1)B_{2,l,e_s}\Delta\rho_{k,2}^{i,s}\right\}.
\end{eqnarray*}
Further, denote 
\[
\rho_{k+1,n}^{i}=\sum_{l=1}^d\rho_{k+1,n}^{i,l},
\]
then we get for $k=j+1,\ldots,J-1$,
\begin{eqnarray*}
\rho_{k+1,2}^{i} & \leq & (1+A_{2}\Delta)\rho_{k,2}^{i}
+2(d-1)C_{1,1}\Delta\rho_{k,2}^{i}+(d-1)^2B_{2}\Delta\rho_{k,2}^{i}.
\end{eqnarray*}
where $A_{2}\doteq\max\limits_{l=1,\ldots,d}A_{2,l}$, $B_{2}\doteq\max\limits_{l,s=1,\ldots,d}B_{2,l,e_s}$ and $C_{1,1}\doteq\max\limits_{l,s=1,\ldots,d}C_{1,1,l,e_l,e_s}$. 
This gives us
\[
\rho_{k+1,2}^{i}\leq (1+\kappa_{1}\Delta)\rho_{k,2}^{i},\quad k=j+1,\ldots,J-1
\]
for some constant $\kappa_{1}>0$,
leading to
\[
\rho_{k,2}^{i}\leq(1+\kappa_{1}\Delta)^{k-j-1}\rho_{j+1,2}^{i}
,\,k=j+1,\ldots,J-1,
\]
where 
\[
\rho_{j+1,2}^{i}=\sum_{s=1}^d\EE\left[\left(\frac{\partial}{\partial y_i}\Phi_{\Delta}^s(x,y+\xi_{j+1})\right)^{2}\right]=\sum_{s=1}^d\sigma_{si}^2(x).
\]
Thus, we obtain the boundedness of 
\[
\frac{\partial}{\partial y_i}u_{\Delta}(t_{j},x,y)=\sum_{s=1}^d\EE\left[\frac{\partial }{\partial x_s}f(G_{J,j}(x,y))\frac{\partial}{\partial y_i}G_{J,j}^s(x,y)\right],
\]
provided that $\sigma_{ki}$ and all the derivatives of order 1 of $f,\mu_k,\sigma_{ki}$ are bounded.

Similar calculations show
that the boundedness of $\sigma_{ki}$ is not necessary to assume
in order to get that $\frac{\partial}{\partial x_l}u_{\Delta}(t_{j},x,y)$ and consequently $g_{j,l}(x)$ for $l\in\left\{1,\ldots,d\right\}$ are bounded (recall~\eqref{eq:13022017a1}).
This yields~(A2) under the assumptions in part~(i)
of Theorem~\ref{func_assump}
(that is, the boundedness of $\sigma_{ki}$ is not needed).

Furthermore, we have, due to $(\sum_{k=1}^da_k)^n\le d^{n-1}\sum_{k=1}^d a_k^n$,
\begin{eqnarray*}
&&\EE\left[\left(\frac{\partial}{\partial y_i}G_{k+1,j}^l(x,y)\right)^{4}\right]\\
&\leq&\EE\left[\left(\frac{\partial}{\partial x_l}\Phi_{\Delta,k+1}^{l}(G_{k,j}(x,y))\frac{\partial}{\partial y_i}G_{k,j}^l(x,y)\right)^4\right.\\
&&\phantom{\EE\left[\right.}+\sum_{s\neq l}\left\{4\left(\frac{\partial}{\partial x_l}\Phi_{\Delta,k+1}^{l}(G_{k,j}(x,y))\frac{\partial}{\partial y_i}G_{k,j}^l(x,y)\right)^3\frac{\partial}{\partial x_s}\Phi_{\Delta,k+1}^{l}(G_{k,j}(x,y))\frac{\partial}{\partial y_i}G_{k,j}^s(x,y)\right.\\
&&\phantom{\EE\left[+\sum_{s\neq l}\left\{\right.\right.}+6(d-1)\left(\frac{\partial}{\partial x_l}\Phi_{\Delta,k+1}^{l}(G_{k,j}(x,y))\frac{\partial}{\partial y_i}G_{k,j}^l(x,y)\frac{\partial}{\partial x_s}\Phi_{\Delta,k+1}^{l}(G_{k,j}(x,y))\frac{\partial}{\partial y_i}G_{k,j}^s(x,y)\right)^2\\
&&\phantom{\EE\left[+\sum_{s\neq l}\left\{\right.\right.}+4(d-1)^2\frac{\partial}{\partial x_l}\Phi_{\Delta,k+1}^{l}(G_{k,j}(x,y))\frac{\partial}{\partial y_i}G_{k,j}^l(x,y)\left(\frac{\partial}{\partial x_s}\Phi_{\Delta,k+1}^{l}(G_{k,j}(x,y))\frac{\partial}{\partial y_i}G_{k,j}^s(x,y)\right)^3\\
&&\phantom{\EE\left[+\sum_{s\neq l}\left\{\right.\right.}\left.\left.+(d-1)^3\left(\frac{\partial}{\partial x_s}\Phi_{\Delta,k+1}^{l}(G_{k,j}(x,y))\frac{\partial}{\partial y_i}G_{k,j}^s(x,y)\right)^4\right\}\right]
\end{eqnarray*}
and thus, due to $4a^3b\le 3a^4+b^4$ and $2a^2b^2\le a^4+b^4$,
\begin{eqnarray*}
\rho_{k+1,4}^{i,l}
&\le & (1+A_{4,l}\Delta)\rho_{k,4}^{i,l}
+\sum_{s\neq l}\left\{C_{3,1,l,e_l,e_s}\Delta(3\rho_{k,4}^{i,l}+\rho_{k,4}^{i,s})+3(d-1)C_{2,2,l,e_l,e_s}\Delta(\rho_{k,4}^{i,l}+\rho_{k,4}^{i,s})\right.\\
&&\phantom{(1+A_{4,l}\Delta)\rho_{k,4}^{i,l}
+\sum_{s\neq l}\left\{\right.}\left.+(d-1)^2C_{1,3,l,e_l,e_s}\Delta(\rho_{k,4}^{i,l}+3\rho_{k,4}^{i,s})+ (d-1)^3B_{4,l,e_s}\Delta\rho_{k,4}^{i,s}\right\}
\end{eqnarray*}
This gives us
\begin{eqnarray*}
\rho_{k+1,4}^{i}
&\le & (1+A_{4}\Delta)\rho_{k,4}^{i}
+4(d-1)C_{3,1}\Delta\rho_{k,4}^{i}+6(d-1)^2C_{2,2}\Delta\rho_{k,4}^{i}\\
&&+4(d-1)^3C_{1,3}\Delta\rho_{k,4}^{i}+ (d-1)^4B_{4}\Delta\rho_{k,4}^{i},
\end{eqnarray*}
where $A_{4}\doteq\max\limits_{l=1,\ldots,d}A_{4,l}$, $B_{4}\doteq\max\limits_{l,s=1,\ldots,d}B_{4,l,e_s}$, $C_{3,1}\doteq\max\limits_{l,s=1,\ldots,d}C_{3,1,l,e_l,e_s}$, $C_{2,2}\doteq\max\limits_{l,s=1,\ldots,d}C_{2,2,l,e_l,e_s}$ and $C_{1,3}\doteq\max\limits_{l,s=1,\ldots,d}C_{1,3,l,e_l,e_s}$.
Hence, we obtain
\[
\rho_{k+1,4}^{i}\leq (1+\kappa_{2}\Delta)\rho_{k,4}^{i},\quad k=j+1,\ldots,J-1
\]
for some constant $\kappa_{2}>0$, leading to 
\[
\rho_{k,4}^{i}\leq(1+\kappa_{2}\Delta)^{k-j-1}\rho_{j+1,4}^{i}
,\,k=j+1,\ldots,J-1,
\]
where 
\[
\rho_{j+1,4}^{i}=\sum_{s=1}^d\EE\left[\left(\frac{\partial}{\partial y_i}\Phi_{\Delta}^s(x,y+\xi_{j+1})\right)^{4}\right]=\sum_{s=1}^d\sigma_{si}^4(x).
\]
Thus, we obtain boundedness of $\rho^i_{k,4}$
uniformly in $x$, $y$, $j$, $k\in\{j+1,\ldots,J\}$ and $J$,
for all $i\in\{1,\ldots,m\}$,
provided $\sigma_{ki}$
and all derivatives of order~1
of $f,\mu_k,\sigma_{ki}$
are bounded.

Now we set\footnote{Notice that thus defined
$\wt G_{J,j}$ is the same
as $G_{J,j}$ of~\eqref{strong:def_G}
(in the proof of Theorem~\ref{th:13022017a1}).}
$\wt G_{J,j}(x):=G_{J,j}(x,0)$ and observe
that similar calculations involving derivatives w.r.t.\ $x_k$
show that the quantities
$$
\EE\left[\left(
\frac\partial{\partial x_k}
\wt G^s_{J,j}(x)
\right)^4\right]
$$
(cf.~with~\eqref{eq:14022017a1})
are all bounded uniformly in
$x$, $J$ and $j\in\{0,\ldots,J-1\}$,
provided all derivatives of order~1
of $f,\mu_k,\sigma_{ki}$
are bounded
(that is, boundedness of $\sigma_{ki}$
is not needed at this point).
Using the identity
$\wt G_{J,j}(X_{\Delta,t_j})=X_{\Delta,T}$
one can check that
\begin{align}
\label{eq:14022017a2}
\cJ_{\wt G_{J,j}}(X_{\Delta,t_j})
=\delta X_{\Delta,T}
\delta X_{\Delta,t_j}^{-1},
\end{align}
where $\cJ_{\wt G_{J,j}}$
denotes the Jacobi matrix of the function $\wt G_{J,j}$.
Recalling the definition
$\zeta_j=(\zeta_{j,1},\ldots,\zeta_{j,d})
:=\nabla f(X_{\Delta,T})\delta X_{\Delta,T}
\delta X_{\Delta,t_j}^{-1}$
of the vector $\zeta_j$,
we get from~\eqref{eq:14022017a2} that
$$
\zeta_{j,k}=\sum_{s=1}^d
\frac\partial{\partial x_s}f(\wt G_{J,j}(X_{\Delta,t_j}))
\frac\partial{\partial x_k}\wt G^s_{J,j}(X_{\Delta,t_j}).
$$
Then we obtain for $k\in\left\{1,\ldots,d\right\}$ and $j\in\left\{1,\ldots,J\right\}$
\begin{align*}
&\Var\left[\zeta_{j,k}\left|\right.X_{\Delta,t_{j-1}}=x\right]\\
&\le \EE\left[\zeta_{j,k}^2\left|\right.X_{\Delta,t_{j-1}}=x\right]\\
&=\EE\left[\left(\sum_{s=1}^d\frac{\partial}{\partial x_s}f(\wt G_{J,j}(X_{\Delta,t_j}))\frac{\partial}{\partial x_k}\wt G_{J,j}^s(X_{\Delta,t_j})\right)^2\left|\right.X_{\Delta,t_{j-1}}=x\right]\\
&\le d\sum_{s=1}^d\EE\left[\left(\frac{\partial}{\partial x_s}f(\wt G_{J,j-1}(x))\frac{\partial}{\partial x_k}\wt G_{J,j}^s(\Phi_{\Delta,j}(x))\right)^2\right]\\
&\le d\sum_{s=1}^d\sqrt{\EE\left[\left(\frac{\partial}{\partial x_s}f(\wt G_{J,j-1}(x))\right)^4\right]\EE\left[\left(\frac{\partial}{\partial x_k}\wt G_{J,j}^s(\Phi_{\Delta_j}(x))\right)^4\right]}.
\end{align*}
Due to the discussion above,
the latter expression is bounded in~$x$,
provided all derivatives of order~1
of $f,\mu_k,\sigma_{ki}$
are bounded.
That is, we get~(A1),
and the proof of part~(i) is completed.

Proceeding with part~(ii), we have
\begin{eqnarray*}
&&\EE\left[\left(\frac{\partial}{\partial y_i}G_{k+1,j}^l(x,y)\right)^{6}\right]\\
&\leq&\EE\left[\left(\frac{\partial}{\partial x_l}\Phi_{\Delta,k+1}^{l}(G_{k,j}(x,y))\frac{\partial}{\partial y_i}G_{k,j}^l(x,y)\right)^6\right.\\
&&\phantom{\EE\left[\right.}+\sum_{s\neq l}\left\{6\left(\frac{\partial}{\partial x_l}\Phi_{\Delta,k+1}^{l}(G_{k,j}(x,y))\frac{\partial}{\partial y_i}G_{k,j}^l(x,y)\right)^5\frac{\partial}{\partial x_s}\Phi_{\Delta,k+1}^{l}(G_{k,j}(x,y))\frac{\partial}{\partial y_i}G_{k,j}^s(x,y)\right.\\
&&\phantom{\EE\left[+\sum_{s\neq l}\left\{\right.\right.}+15(d-1)\left(\frac{\partial}{\partial x_l}\Phi_{\Delta,k+1}^{l}(G_{k,j}(x,y))\frac{\partial}{\partial y_i}G_{k,j}^l(x,y)\right)^4\left(\frac{\partial}{\partial x_s}\Phi_{\Delta,k+1}^{l}(G_{k,j}(x,y))\frac{\partial}{\partial y_i}G_{k,j}^s(x,y)\right)^2\\
&&\phantom{\EE\left[+\sum_{s\neq l}\left\{\right.\right.}+20(d-1)^2\left(\frac{\partial}{\partial x_l}\Phi_{\Delta,k+1}^{l}(G_{k,j}(x,y))\frac{\partial}{\partial y_i}G_{k,j}^l(x,y)\frac{\partial}{\partial x_s}\Phi_{\Delta,k+1}^{l}(G_{k,j}(x,y))\frac{\partial}{\partial y_i}G_{k,j}^s(x,y)\right)^3\\
&&\phantom{\EE\left[+\sum_{s\neq l}\left\{\right.\right.}+15(d-1)^3\left(\frac{\partial}{\partial x_l}\Phi_{\Delta,k+1}^{l}(G_{k,j}(x,y))\frac{\partial}{\partial y_i}G_{k,j}^l(x,y)\right)^2\left(\frac{\partial}{\partial x_s}\Phi_{\Delta,k+1}^{l}(G_{k,j}(x,y))\frac{\partial}{\partial y_i}G_{k,j}^s(x,y)\right)^4\\
&&\phantom{\EE\left[+\sum_{s\neq l}\left\{\right.\right.}+6(d-1)^4\frac{\partial}{\partial x_l}\Phi_{\Delta,k+1}^{l}(G_{k,j}(x,y))\frac{\partial}{\partial y_i}G_{k,j}^l(x,y)\left(\frac{\partial}{\partial x_s}\Phi_{\Delta,k+1}^{l}(G_{k,j}(x,y))\frac{\partial}{\partial y_i}G_{k,j}^s(x,y)\right)^5\\
&&\phantom{\EE\left[+\sum_{s\neq l}\left\{\right.\right.}\left.\left.+(d-1)^5\left(\frac{\partial}{\partial x_s}\Phi_{\Delta,k+1}^{l}(G_{k,j}(x,y))\frac{\partial}{\partial y_i}G_{k,j}^s(x,y)\right)^6\right\}\right]
\end{eqnarray*}
and thus, due to $6a^5b\le 5a^6+b^6$, $3a^4b^2\le 2a^6+b^6$ and $2a^3b^3\le a^6+b^6$,
\begin{eqnarray*}
\rho_{k+1,6}^{i,l}
&\le & (1+A_{6,l}\Delta)\rho_{k,6}^{i,l}
+\sum_{s\neq l}\left\{C_{5,1,l,e_l,e_s}\Delta(5\rho_{k,6}^{i,l}+\rho_{k,6}^{i,s})+5(d-1)C_{4,2,l,e_l,e_s}\Delta(2\rho_{k,6}^{i,l}+\rho_{k,6}^{i,s})\right.\\
&&\phantom{(1+A_{4,l}\Delta)\rho_{k,6}^{i,l}
+\sum_{s\neq l}\left\{\right.}+10(d-1)^2C_{3,3,l,e_l,e_s}\Delta(\rho_{k,6}^{i,l}+\rho_{k,6}^{i,s})+5(d-1)^3C_{2,4,l,e_l,e_s}\Delta(\rho_{k,6}^{i,l}+2\rho_{k,6}^{i,s})\\
&&\phantom{(1+A_{4,l}\Delta)\rho_{k,6}^{i,l}
+\sum_{s\neq l}\left\{\right.}\left.+(d-1)^4C_{1,5,l,e_l,e_s}\Delta(\rho_{k,6}^{i,l}+5\rho_{k,6}^{i,s})+ (d-1)^5B_{6,l,e_s}\Delta\rho_{k,6}^{i,s}\right\}.
\end{eqnarray*}
This gives us
\begin{eqnarray*}
\rho_{k+1,6}^{i}
&\le & (1+A_{6}\Delta)\rho_{k,6}^{i}
+6(d-1)C_{5,1}\Delta\rho_{k,6}^{i}+15(d-1)^2C_{4,2}\Delta\rho_{k,6}^{i}+20(d-1)^3C_{3,3}\Delta\rho_{k,6}^{i}\\
&&+15(d-1)^4C_{2,4}\Delta\rho_{k,6}^{i}+6(d-1)^5C_{1,5}\Delta\rho_{k,6}^{i}+ (d-1)^6B_{6}\Delta\rho_{k,6}^{i},
\end{eqnarray*}
where $A_{6}\doteq\max\limits_{l=1,\ldots,d}A_{6,l}$, $B_{6}\doteq\max\limits_{l,s=1,\ldots,d}B_{6,l,e_s}$, $C_{5,1}\doteq\max\limits_{l,s=1,\ldots,d}C_{5,1,l,e_l,e_s}$, $C_{4,2}\doteq\max\limits_{l,s=1,\ldots,d}C_{4,2,l,e_l,e_s}$, $C_{3,3}\doteq\max\limits_{l,s=1,\ldots,d}C_{3,3,l,e_l,e_s}$, $C_{2,4}\doteq\max\limits_{l,s=1,\ldots,d}C_{2,4,l,e_l,e_s}$ and $C_{1,5}\doteq\max\limits_{l,s=1,\ldots,d}C_{1,5,l,e_l,e_s}$.
Hence, we obtain
\[
\rho_{k+1,6}^{i}\leq (1+\kappa_{3}\Delta)\rho_{k,6}^{i},\quad k=j+1,\ldots,J-1
\]
for some constant $\kappa_{3}>0$, leading to 
\[
\rho_{k,6}^{i}\leq(1+\kappa_{3}\Delta)^{k-j-1}\rho_{j+1,6}^{i}
,\,k=j+1,\ldots,J-1,
\]
where 
\[
\rho_{j+1,6}^{i}=\sum_{s=1}^d\EE\left[\left(\frac{\partial}{\partial y_i}\Phi_{\Delta}^s(x,y+\xi_{j+1})\right)^{6}\right]=\sum_{s=1}^d\sigma_{si}^6(x).
\]

Moreover, we have
\begin{eqnarray*}
&&\EE\left[\left(\frac{\partial}{\partial y_i}G_{k+1,j}^l(x,y)\right)^{8}\right]\\
&\leq&\EE\left[\left(\frac{\partial}{\partial x_l}\Phi_{\Delta,k+1}^{l}(G_{k,j}(x,y))\frac{\partial}{\partial y_i}G_{k,j}^l(x,y)\right)^8\right.\\
&&\phantom{\EE\left[\right.}+\sum_{s\neq l}\left\{8\left(\frac{\partial}{\partial x_l}\Phi_{\Delta,k+1}^{l}(G_{k,j}(x,y))\frac{\partial}{\partial y_i}G_{k,j}^l(x,y)\right)^7\frac{\partial}{\partial x_s}\Phi_{\Delta,k+1}^{l}(G_{k,j}(x,y))\frac{\partial}{\partial y_i}G_{k,j}^s(x,y)\right.\\
&&\phantom{\EE\left[+\sum_{s\neq l}\left\{\right.\right.}+28(d-1)\left(\frac{\partial}{\partial x_l}\Phi_{\Delta,k+1}^{l}(G_{k,j}(x,y))\frac{\partial}{\partial y_i}G_{k,j}^l(x,y)\right)^6\left(\frac{\partial}{\partial x_s}\Phi_{\Delta,k+1}^{l}(G_{k,j}(x,y))\frac{\partial}{\partial y_i}G_{k,j}^s(x,y)\right)^2\\
&&\phantom{\EE\left[+\sum_{s\neq l}\left\{\right.\right.}+56(d-1)^2\left(\frac{\partial}{\partial x_l}\Phi_{\Delta,k+1}^{l}(G_{k,j}(x,y))\frac{\partial}{\partial y_i}G_{k,j}^l(x,y)\right)^5\left(\frac{\partial}{\partial x_s}\Phi_{\Delta,k+1}^{l}(G_{k,j}(x,y))\frac{\partial}{\partial y_i}G_{k,j}^s(x,y)\right)^3\\
&&\phantom{\EE\left[+\sum_{s\neq l}\left\{\right.\right.}+70(d-1)^3\left(\frac{\partial}{\partial x_l}\Phi_{\Delta,k+1}^{l}(G_{k,j}(x,y))\frac{\partial}{\partial y_i}G_{k,j}^l(x,y)\frac{\partial}{\partial x_s}\Phi_{\Delta,k+1}^{l}(G_{k,j}(x,y))\frac{\partial}{\partial y_i}G_{k,j}^s(x,y)\right)^4\\
&&\phantom{\EE\left[+\sum_{s\neq l}\left\{\right.\right.}+56(d-1)^4\left(\frac{\partial}{\partial x_l}\Phi_{\Delta,k+1}^{l}(G_{k,j}(x,y))\frac{\partial}{\partial y_i}G_{k,j}^l(x,y)\right)^3\left(\frac{\partial}{\partial x_s}\Phi_{\Delta,k+1}^{l}(G_{k,j}(x,y))\frac{\partial}{\partial y_i}G_{k,j}^s(x,y)\right)^5\\
&&\phantom{\EE\left[+\sum_{s\neq l}\left\{\right.\right.}+28(d-1)^5\left(\frac{\partial}{\partial x_l}\Phi_{\Delta,k+1}^{l}(G_{k,j}(x,y))\frac{\partial}{\partial y_i}G_{k,j}^l(x,y)\right)^2\left(\frac{\partial}{\partial x_s}\Phi_{\Delta,k+1}^{l}(G_{k,j}(x,y))\frac{\partial}{\partial y_i}G_{k,j}^s(x,y)\right)^6\\
&&\phantom{\EE\left[+\sum_{s\neq l}\left\{\right.\right.}+8(d-1)^6\frac{\partial}{\partial x_l}\Phi_{\Delta,k+1}^{l}(G_{k,j}(x,y))\frac{\partial}{\partial y_i}G_{k,j}^l(x,y)\left(\frac{\partial}{\partial x_s}\Phi_{\Delta,k+1}^{l}(G_{k,j}(x,y))\frac{\partial}{\partial y_i}G_{k,j}^s(x,y)\right)^7\\
&&\phantom{\EE\left[+\sum_{s\neq l}\left\{\right.\right.}\left.\left.+(d-1)^7\left(\frac{\partial}{\partial x_s}\Phi_{\Delta,k+1}^{l}(G_{k,j}(x,y))\frac{\partial}{\partial y_i}G_{k,j}^s(x,y)\right)^8\right\}\right]
\end{eqnarray*}
and thus, due to $8a^7b\le 7a^8+b^8$, $4a^6b^2\le 3a^8+b^8$, $8a^5b^3\le 5a^8+3b^8$ and $2a^4b^4\le a^8+b^8$,
\begin{eqnarray*}
\rho_{k+1,8}^{i,l}
&\le & (1+A_{8,l}\Delta)\rho_{k,8}^{i,l}
+\sum_{s\neq l}\left\{C_{7,1,l,e_l,e_s}\Delta(7\rho_{k,8}^{i,l}+\rho_{k,8}^{i,s})+7(d-1)C_{6,2,l,e_l,e_s}\Delta(3\rho_{k,8}^{i,l}+\rho_{k,8}^{i,s})\right.\\
&&\phantom{(1+A_{4,l}\Delta)\rho_{k,8}^{i,l}
+\sum_{s\neq l}\left\{\right.}+7(d-1)^2C_{5,3,l,e_l,e_s}\Delta(5\rho_{k,8}^{i,l}+3\rho_{k,8}^{i,s})+35(d-1)^3C_{4,4,l,e_l,e_s}\Delta(\rho_{k,8}^{i,l}+\rho_{k,8}^{i,s})\\
&&\phantom{(1+A_{4,l}\Delta)\rho_{k,8}^{i,l}
+\sum_{s\neq l}\left\{\right.}+7(d-1)^4C_{3,5,l,e_l,e_s}\Delta(3\rho_{k,8}^{i,l}+5\rho_{k,8}^{i,s})+7(d-1)^5C_{2,6,l,e_l,e_s}\Delta(\rho_{k,8}^{i,l}+3\rho_{k,8}^{i,s})\\
&&\phantom{(1+A_{4,l}\Delta)\rho_{k,8}^{i,l}
+\sum_{s\neq l}\left\{\right.}\left.+(d-1)^6C_{1,7,l,e_l,e_s}\Delta(\rho_{k,8}^{i,l}+7\rho_{k,8}^{i,s})+ (d-1)^7B_{8,l,e_s}\Delta\rho_{k,8}^{i,s}\right\}.
\end{eqnarray*}
This gives us
\begin{eqnarray*}
\rho_{k+1,8}^{i}
&\le & (1+A_{8}\Delta)\rho_{k,8}^{i}
+8(d-1)C_{7,1}\Delta\rho_{k,8}^{i}+28(d-1)^2C_{6,2}\Delta\rho_{k,8}^{i}+56(d-1)^3C_{5,3}\Delta\rho_{k,8}^{i}\\
&&+70(d-1)^4C_{4,4}\Delta\rho_{k,8}^{i}+56(d-1)^5C_{3,5}\Delta\rho_{k,8}^{i}+28(d-1)^6C_{2,6}\Delta\rho_{k,8}^{i}\\
&&+8(d-1)^7C_{1,7}\Delta\rho_{k,8}^{i}+ (d-1)^8B_{8}\Delta\rho_{k,8}^{i},
\end{eqnarray*}
where $A_{8}\doteq\max\limits_{l=1,\ldots,d}A_{8,l}$, $B_{8}\doteq\max\limits_{l,s=1,\ldots,d}B_{8,l,e_s}$, $C_{7,1}\doteq\max\limits_{l,s=1,\ldots,d}C_{7,1,l,e_l,e_s}$, $C_{6,2}\doteq\max\limits_{l,s=1,\ldots,d}C_{6,2,l,e_l,e_s}$, $C_{5,3}\doteq\max\limits_{l,s=1,\ldots,d}C_{5,3,l,e_l,e_s}$, $C_{4,4}\doteq\max\limits_{l,s=1,\ldots,d}C_{4,4,l,e_l,e_s}$, $C_{3,5}\doteq\max\limits_{l,s=1,\ldots,d}C_{3,5,l,e_l,e_s}$, $C_{2,6}\doteq\max\limits_{l,s=1,\ldots,d}C_{2,6,l,e_l,e_s}$ and $C_{1,7}\doteq\max\limits_{l,s=1,\ldots,d}C_{1,7,l,e_l,e_s}$.
Hence, we obtain
\[
\rho_{k+1,8}^{i}\leq (1+\kappa_{4}\Delta)\rho_{k,8}^{i},\quad k=j+1,\ldots,J-1
\]
for some constant $\kappa_{4}>0$, leading to 
\[
\rho_{k,8}^{i}\leq(1+\kappa_{4}\Delta)^{k-j-1}\rho_{j+1,8}^{i}
,\,k=j+1,\ldots,J-1,
\]
where 
\[
\rho_{j+1,8}^{i}=\sum_{s=1}^d\EE\left[\left(\frac{\partial}{\partial y_i}\Phi_{\Delta}^s(x,y+\xi_{j+1})\right)^{8}\right]=\sum_{s=1}^d\sigma_{si}^8(x).
\]

Next, we have for some $i,o\in\left\{1,\ldots,m\right\}$ and $l\in\left\{1,\ldots,d\right\}$
\begin{align*}
\frac{\partial^2}{\partial y_i\partial y_o}G_{k+1,j}^l(x,y)=&\sum_{s=1}^d\frac{\partial}{\partial x_s}\Phi_{\Delta,k+1}^{l}(G_{k,j}(x,y))\frac{\partial^2}{\partial y_i\partial y_o}G_{k,j}^s(x,y)\\
&+\sum_{s,u=1}^d\frac{\partial^2}{\partial x_s\partial x_u}\Phi_{\Delta,k+1}^{l}(G_{k,j}(x,y))\frac{\partial}{\partial y_i}G_{k,j}^s(x,y)\frac{\partial}{\partial y_o}G_{k,j}^u(x,y)
\end{align*}
and $\frac{\partial^2}{\partial y_i\partial y_o}G_{j+1,j}^s(x,y)=\frac{\partial^2}{\partial y_i\partial y_o}\Phi_{\Delta}^s(x,y+\xi_{j+1})$.
Hence 
\begin{eqnarray*}
&&\EE\left[\left(\frac{\partial^2}{\partial y_i\partial y_o}G_{k+1,j}^l(x,y)\right)^{2}\right]\\
&\leq&\EE\left[\left(\frac{\partial}{\partial x_l}\Phi_{\Delta,k+1}^{l}(G_{k,j}(x,y))\frac{\partial^2}{\partial y_i\partial y_o}G_{k,j}^l(x,y)\right)^2\right.\\
&&\phantom{\EE\left[\right.}+\sum_{s\neq l}\left\{2\frac{\partial}{\partial x_l}\Phi_{\Delta,k+1}^{l}(G_{k,j}(x,y))\frac{\partial}{\partial x_s}\Phi_{\Delta,k+1}^{l}(G_{k,j}(x,y))\frac{\partial^2}{\partial y_i\partial y_o}G_{k,j}^l(x,y)\frac{\partial^2}{\partial y_i\partial y_o}G_{k,j}^s(x,y)\right.\\
&&\phantom{\EE\left[+\sum_{s\neq l}\left\{\right.\right.}\left.+(d-1)\left(\frac{\partial}{\partial x_s}\Phi_{\Delta,k+1}^{l}(G_{k,j}(x,y))\frac{\partial^2}{\partial y_i\partial y_o}G_{k,j}^s(x,y)\right)^2\right\}\\
&&\phantom{\EE\left[\right.}+2\sum_{s,u,v=1}^d\frac{\partial}{\partial x_v}\Phi_{\Delta,k+1}^{l}(G_{k,j}(x,y))\frac{\partial^2}{\partial x_s\partial x_u}\Phi_{\Delta,k+1}^{l}(G_{k,j}(x,y))\frac{\partial^2}{\partial y_i\partial y_o}G_{k,j}^v(x,y)\frac{\partial}{\partial y_i}G_{k,j}^s(x,y)\frac{\partial}{\partial y_o}G_{k,j}^u(x,y)\\
&&\phantom{\EE\left[\right.}\left.+d^2\sum_{s,u=1}^d\left(\frac{\partial^2}{\partial x_s\partial x_u}\Phi_{\Delta,k+1}^{l}(G_{k,j}(x,y))\frac{\partial}{\partial y_i}G_{k,j}^s(x,y)\frac{\partial}{\partial y_o}G_{k,j}^u(x,y)\right)^2\right]
\end{eqnarray*}
Denote 
\[
\psi_{k+1,n}^{i,o,s}=\EE\left[\left(\frac{\partial^2}{\partial y_i\partial y_o}G_{k+1,j}^s(x,y)\right)^{n}\right]
\]
and $e_{s,u}\doteq e_s+e_u$, then we get, due to 
$$2\EE\left[XYZ\right]\le 2\sqrt{\EE\left[X^2\right]}\sqrt[4]{\EE\left[Y^4\right]}\sqrt[4]{\EE\left[Z^4\right]}\le \EE\left[X^2\right]+\sqrt{\EE\left[Y^4\right]}\sqrt{\EE\left[Z^4\right]}\le \EE\left[X^2\right]+\frac{1}{2}\left(\EE\left[Y^4\right]+\EE\left[Z^4\right]\right),
$$
for $k=j+1,\ldots,J-1$,
\begin{eqnarray*}
\psi_{k+1,2}^{i,o,l} & \leq & (1+A_{2,l}\Delta)\psi_{k,2}^{i,o,l}
+\sum_{s\neq l}\left\{C_{1,1,l,e_l,e_s}\Delta(\psi_{k,2}^{i,o,l}+\psi_{k,2}^{i,o,s})+(d-1)B_{2,l,e_s}\Delta\psi_{k,2}^{i,o,s}\right\}\\
&&+\sum_{s,u,v=1}^d C_{1,1,l,e_v,e_{s,u}}\Delta\left(\psi_{k,2}^{i,o,v}+\frac{1}{2}\left(\rho_{k,4}^{i,s}+\rho_{k,4}^{o,u}\right)\right)+d^2\sum_{s,u=1}^dB_{2,l,e_{s,u}}\Delta\frac{1}{2}\left(\rho_{k,4}^{i,s}+\rho_{k,4}^{o,u}\right).
\end{eqnarray*}
Further, denote 
\[
\psi_{k+1,n}^{i,o}=\sum_{l=1}^d\psi_{k+1,n}^{i,o,l},
\]
then we get for $k=j+1,\ldots,J-1$,
\begin{eqnarray*}
\psi_{k+1,2}^{i,o} & \leq & (1+A_{2}\Delta)\psi_{k,2}^{i,o}
+2(d-1)C_{1,1}\Delta\psi_{k,2}^{i,o}+(d-1)^2B_{2}\Delta\psi_{k,2}^{i,o}\\
&&+d^3\tilde{C}_{1,1}\Delta\left(\psi_{k,2}^{i,o}+\frac{1}{2}\left(\rho_{k,4}^{i}+\rho_{k,4}^{o}\right)\right)+d^4\tilde{B}_2\Delta\frac{1}{2}\left(\rho_{k,4}^{i}+\rho_{k,4}^{o}\right).
\end{eqnarray*}
where $\tilde{C}_{1,1}\doteq\max\limits_{l,s,u,v=1,\ldots,d}C_{1,1,l,e_v,e_{s,u}}$ and $\tilde{B}_{2}\doteq\max\limits_{l,s,u=1,\ldots,d}B_{2,l,e_{s,u}}$. 
This gives us
\[
\psi_{k+1,2}^{i,o}\leq (1+\kappa_{5}\Delta)\psi_{k,2}^{i,o}+\kappa_6,\quad k=j+1,\ldots,J-1
\]
for some constants $\kappa_{5},\kappa_6>0$,
leading to
\[
\psi_{k,2}^{i,o}\leq(1+\kappa_{5}\Delta)^{k-j-1}\psi_{j+1,2}^{i,o}+\kappa_7=\kappa_7
,\,k=j+1,\ldots,J-1,
\]
where $\kappa_7>0$ and
\[
\psi_{j+1,2}^{i,o}=\sum_{s=1}^d\EE\left[\left(\frac{\partial^2}{\partial y_i\partial y_o}\Phi_{\Delta}^s(x,y+\xi_{j+1})\right)^{2}\right]=0.
\]
Thus, we obtain the boundedness of 
\begin{eqnarray*}
\frac{\partial^2}{\partial y_i\partial y_o}u_{\Delta}(t_{j},x,y)&=&\EE\left[\sum_{s=1}^d\frac{\partial }{\partial x_s}f(G_{J,j}(x,y))\frac{\partial^2}{\partial y_i\partial y_o}G_{J,j}^s(x,y)\right.\\
&&\phantom{\EE\left[\right.}\left.+\sum_{s,u=1}^d\frac{\partial^2 }{\partial x_s\partial x_u}f(G_{J,j}(x,y))\frac{\partial}{\partial y_i}G_{J,j}^s(x,y)\frac{\partial}{\partial y_o}G_{J,j}^u(x,y)\right],
\end{eqnarray*}
provided that $\sigma_{ki}$ and all the derivatives of order 1 and 2 for $f,\mu_k,\sigma_{ki}$ are bounded.

Moreover, we have
\begin{eqnarray*}
&&\EE\left[\left(\frac{\partial^2}{\partial y_i\partial y_o}G_{k+1,j}^l(x,y)\right)^{4}\right]\\
&\leq&\EE\left[\left(\frac{\partial}{\partial x_l}\Phi_{\Delta,k+1}^{l}(G_{k,j}(x,y))\frac{\partial^2}{\partial y_i\partial y_o}G_{k,j}^l(x,y)\right)^4\right.\\
&&\phantom{\EE\left[\right.}+\sum_{s\neq l}\left\{4\left(\frac{\partial}{\partial x_l}\Phi_{\Delta,k+1}^{l}(G_{k,j}(x,y))\frac{\partial^2}{\partial y_i\partial y_o}G_{k,j}^l(x,y)\right)^3\frac{\partial}{\partial x_s}\Phi_{\Delta,k+1}^{l}(G_{k,j}(x,y))\frac{\partial^2}{\partial y_i\partial y_o}G_{k,j}^s(x,y)\right.\\
&&\phantom{\EE\left[+\sum_{s\neq l}\left\{\right.\right.}+6(d-1)\left(\frac{\partial}{\partial x_l}\Phi_{\Delta,k+1}^{l}(G_{k,j}(x,y))\frac{\partial^2}{\partial y_i\partial y_o}G_{k,j}^l(x,y)\frac{\partial}{\partial x_s}\Phi_{\Delta,k+1}^{l}(G_{k,j}(x,y))\frac{\partial^2}{\partial y_i\partial y_o}G_{k,j}^s(x,y)\right)^2\\
&&\phantom{\EE\left[+\sum_{s\neq l}\left\{\right.\right.}+4(d-1)^2\frac{\partial}{\partial x_l}\Phi_{\Delta,k+1}^{l}(G_{k,j}(x,y))\frac{\partial^2}{\partial y_i\partial y_o}G_{k,j}^l(x,y)\left(\frac{\partial}{\partial x_s}\Phi_{\Delta,k+1}^{l}(G_{k,j}(x,y))\frac{\partial^2}{\partial y_i\partial y_o}G_{k,j}^s(x,y)\right)^3\\
&&\phantom{\EE\left[+\sum_{s\neq l}\left\{\right.\right.}\left.+(d-1)^3\left(\frac{\partial}{\partial x_s}\Phi_{\Delta,k+1}^{l}(G_{k,j}(x,y))\frac{\partial^2}{\partial y_i\partial y_o}G_{k,j}^s(x,y)\right)^4\right\}\\
&&\phantom{\EE\left[\right.}+\sum_{s,u,v=1}^d\left\{4d^2\left(\frac{\partial}{\partial x_v}\Phi_{\Delta,k+1}^{l}(G_{k,j}(x,y))\frac{\partial^2}{\partial y_i\partial y_o}G_{k,j}^v(x,y)\right)^3\frac{\partial^2}{\partial x_s\partial x_u}\Phi_{\Delta,k+1}^{l}(G_{k,j}(x,y))\right.\\
&&\phantom{\EE\left[+\sum_{s,u,v=1}^d\left\{\right.\right.}\cdot\frac{\partial}{\partial y_i}G_{k,j}^s(x,y)\frac{\partial}{\partial y_o}G_{k,j}^u(x,y)\\
&&\phantom{\EE\left[+\sum_{s,u,v=1}^d\left\{\right.\right.}+6d^3\left(\frac{\partial}{\partial x_v}\Phi_{\Delta,k+1}^{l}(G_{k,j}(x,y))\frac{\partial^2}{\partial y_i\partial y_o}G_{k,j}^v(x,y)\frac{\partial^2}{\partial x_s\partial x_u}\Phi_{\Delta,k+1}^{l}(G_{k,j}(x,y))\right)^2\\
&&\phantom{\EE\left[+\sum_{s,u,v=1}^d\left\{+\right.\right.}\cdot \left(\frac{\partial}{\partial y_i}G_{k,j}^s(x,y)\frac{\partial}{\partial y_o}G_{k,j}^u(x,y)\right)^2\\
&&\phantom{\EE\left[+\sum_{s,u,v=1}^d\left\{\right.\right.}+4d^4\frac{\partial}{\partial x_v}\Phi_{\Delta,k+1}^{l}(G_{k,j}(x,y))\frac{\partial^2}{\partial y_i\partial y_o}G_{k,j}^v(x,y)\\
&&\phantom{\EE\left[+\sum_{s,u,v=1}^d\left\{+\right.\right.}\left.\cdot\left(\frac{\partial^2}{\partial x_s\partial x_u}\Phi_{\Delta,k+1}^{l}(G_{k,j}(x,y))\frac{\partial}{\partial y_i}G_{k,j}^s(x,y)\frac{\partial}{\partial y_o}G_{k,j}^u(x,y)\right)^3\right\}\\
&&\phantom{\EE\left[\right.}\left.+d^6\sum_{s,u}^d\left(\frac{\partial^2}{\partial x_s\partial x_u}\Phi_{\Delta,k+1}^{l}(G_{k,j}(x,y))\frac{\partial}{\partial y_i}G_{k,j}^s(x,y)\frac{\partial}{\partial y_o}G_{k,j}^u(x,y)\right)^4\right]
\end{eqnarray*}
and thus, due to $4a^3bc\le3a^4+\frac{1}{2}\left(b^8+c^8\right)$, $2a^2b^2c^2\le a^4+\frac{1}{2}\left(b^8+c^8\right)$ and $4ab^3c^3\le a^4+\frac{3}{2}\left(b^8+c^8\right)$,
\begin{eqnarray*}
\psi_{k+1,4}^{i,o,l}
&\le & (1+A_{4,l}\Delta)\psi_{k,4}^{i,o,l}
+\sum_{s\neq l}\left\{C_{3,1,l,e_l,e_s}\Delta(3\psi_{k,4}^{i,o,l}+\psi_{k,4}^{i,o,s})+3(d-1)C_{2,2,l,e_l,e_s}\Delta(\psi_{k,4}^{i,o,l}+\psi_{k,4}^{i,o,s})\right.\\
&&\phantom{(1+A_{4,l}\Delta)\rho_{k,4}^{i,l}
+\sum_{s\neq l}\left\{\right.}\left.+(d-1)^2C_{1,3,l,e_l,e_s}\Delta(\psi_{k,4}^{i,o,l}+3\psi_{k,4}^{i,o,s})+ (d-1)^3B_{4,l,e_s}\Delta\psi_{k,4}^{i,o,s}\right\}\\
&&+\sum_{s,u,v=1}^d\left\{d^2C_{3,1,l,e_v,e_{s,u}}\Delta\left(3\psi_{k,4}^{i,o,v}+\frac{1}{2}\left(\rho_{k,8}^{i,s}+\rho_{k,8}^{o,u}\right)\right)+3d^3C_{2,2,l,e_v,e_{s,u}}\Delta\left(\psi_{k,4}^{i,o,v}+\frac{1}{2}\left(\rho_{k,8}^{i,s}+\rho_{k,8}^{o,u}\right)\right)\right.\\
&&\phantom{+\sum_{s,u,v=1}^d\left\{\right.}\left.+d^4C_{1,3,l,e_v,e_{s,u}}\Delta\left(\psi_{k,4}^{i,o,v}+\frac{3}{2}\left(\rho_{k,8}^{i,s}+\rho_{k,8}^{o,u}\right)\right)\right\}+d^6\sum_{s,u=1}^dB_{4,l,e_{s,u}}\Delta\frac{1}{2}\left(\rho_{k,8}^{i,s}+\rho_{k,8}^{o,u}\right).
\end{eqnarray*}
This gives us
\begin{eqnarray*}
\psi_{k+1,4}^{i,o}
&\le & (1+A_{4}\Delta)\psi_{k,4}^{i,o}
+4(d-1)C_{3,1}\Delta\psi_{k,4}^{i,o}+6(d-1)^2C_{2,2}\Delta\psi_{k,4}^{i,o}\\
&&+4(d-1)^3C_{1,3}\Delta\rho_{k,4}^{i}+ (d-1)^4B_{4}\Delta\psi_{k,4}^{i,o}+d^5\tilde{C}_{3,1}\Delta\left(3\psi_{k,4}^{i,o}+\frac{1}{2}\left(\rho_{k,8}^{i}+\rho_{k,8}^{o}\right)\right)\\
&&+3d^6\tilde{C}_{2,2}\Delta\left(\psi_{k,4}^{i,o}+\frac{1}{2}\left(\rho_{k,8}^{i}+\rho_{k,8}^{o}\right)\right)+d^7\tilde{C}_{1,3}\Delta\left(\psi_{k,4}^{i,o}+\frac{3}{2}\left(\rho_{k,8}^{i}+\rho_{k,8}^{o}\right)\right)+d^8\tilde{B}_{4}\Delta\frac{1}{2}\left(\rho_{k,8}^{i}+\rho_{k,8}^{o}\right),
\end{eqnarray*}
where $\tilde C_{3,1}\doteq\max\limits_{l,s,u,v=1,\ldots,d}C_{3,1,l,e_v,e_{s,u}}$, $\tilde C_{2,2}\doteq\max\limits_{l,s,u,v=1,\ldots,d}C_{2,2,l,e_v,e_{s,u}}$, $\tilde C_{1,3}\doteq\max\limits_{l,s,u,v=1,\ldots,d}C_{1,3,l,e_v,e_{s,u}}$ and $\tilde{B}_{4}\doteq\max\limits_{l,s,u=1,\ldots,d}B_{4,l,e_{s,u}}$.
Hence, we obtain
\[
\psi_{k+1,4}^{i,o}\leq (1+\kappa_{8}\Delta)\psi_{k,4}^{i,o}+\kappa_9,\quad k=j+1,\ldots,J-1
\]
for some constants $\kappa_{8},\kappa_9>0$, leading to 
\[
\psi_{k,4}^{i,o}\leq(1+\kappa_{8}\Delta)^{k-j-1}\psi_{j+1,4}^{i,o}+\kappa_{10}=\kappa_{10}
,\,k=j+1,\ldots,J-1,
\]
where $\kappa_{10}>0$ and
\[
\psi_{j+1,4}^{i,o}=\sum_{s=1}^d\EE\left[\left(\frac{\partial^2}{\partial y_i\partial y_o}\Phi_{\Delta}^s(x,y+\xi_{j+1})\right)^{4}\right]=0.
\]

Next, we have for some $i,o,r\in\left\{1,\ldots,m\right\}$ and $l\in\left\{1,\ldots,d\right\}$
\begin{align*}
&\frac{\partial^3}{\partial y_i\partial y_o\partial y_r}G_{k+1,j}^l(x,y)\\
=&\sum_{s=1}^d\frac{\partial}{\partial x_s}\Phi_{\Delta,k+1}^{l}(G_{k,j}(x,y))\frac{\partial^3}{\partial y_i\partial y_o\partial y_r}G_{k,j}^s(x,y)\\
&+\sum_{s,u=1}^d\frac{\partial^2}{\partial x_s\partial x_u}\Phi_{\Delta,k+1}^{l}(G_{k,j}(x,y))\left(\frac{\partial^2}{\partial y_i\partial y_o}G_{k,j}^s(x,y)\frac{\partial}{\partial y_r}G_{k,j}^u(x,y)+\frac{\partial^2}{\partial y_i\partial y_r}G_{k,j}^s(x,y)\frac{\partial}{\partial y_o}G_{k,j}^u(x,y)\right.\\
&\phantom{+\sum_{s,u=1}^d\frac{\partial^2}{\partial x_s\partial x_u}\Phi_{\Delta,k+1}^{l}(G_{k,j}(x,y))\left(\right.}\left.+\frac{\partial}{\partial y_i}G_{k,j}^s(x,y)\frac{\partial^2}{\partial y_o\partial y_r}G_{k,j}^u(x,y)\right)\\
&+\sum_{s,u,v=1}^d\frac{\partial^3}{\partial x_s\partial x_u\partial x_v}\Phi_{\Delta,k+1}^{l}(G_{k,j}(x,y))\frac{\partial}{\partial y_i}G_{k,j}^s(x,y)\frac{\partial}{\partial y_o}G_{k,j}^u(x,y)\frac{\partial}{\partial y_r}G_{k,j}^v(x,y)
\end{align*}
and $\frac{\partial^3}{\partial y_i\partial y_o\partial y_r}G_{j+1,j}^s(x,y)=\frac{\partial^3}{\partial y_i\partial y_o\partial y_r}\Phi_{\Delta}^s(x,y+\xi_{j+1})$.
Hence 
\begin{eqnarray*}
&&\EE\left[\left(\frac{\partial^3}{\partial y_i\partial y_o\partial y_r}G_{k+1,j}^l(x,y)\right)^{2}\right]\\
&\leq&\EE\left[\left(\frac{\partial}{\partial x_l}\Phi_{\Delta,k+1}^{l}(G_{k,j}(x,y))\frac{\partial^3}{\partial y_i\partial y_o\partial y_r}G_{k,j}^l(x,y)\right)^2\right.\\
&&\phantom{\EE\left[\right.}+\sum_{s\neq l}\left\{2\frac{\partial}{\partial x_l}\Phi_{\Delta,k+1}^{l}(G_{k,j}(x,y))\frac{\partial}{\partial x_s}\Phi_{\Delta,k+1}^{l}(G_{k,j}(x,y))\frac{\partial^3}{\partial y_i\partial y_o\partial y_r}G_{k,j}^l(x,y)\frac{\partial^3}{\partial y_i\partial y_o\partial y_r}G_{k,j}^s(x,y)\right.\\
&&\phantom{\EE\left[+\sum_{s\neq l}\left\{\right.\right.}\left.+(d-1)\left(\frac{\partial}{\partial x_s}\Phi_{\Delta,k+1}^{l}(G_{k,j}(x,y))\frac{\partial^3}{\partial y_i\partial y_o\partial y_r}G_{k,j}^s(x,y)\right)^2\right\}\\
&&\phantom{\EE\left[\right.}+2\sum_{s,u,v=1}^d\frac{\partial^2}{\partial x_s\partial x_u}\Phi_{\Delta,k+1}^{l}(G_{k,j}(x,y))\frac{\partial}{\partial x_v}\Phi_{\Delta,k+1}^{l}(G_{k,j}(x,y))\frac{\partial^3}{\partial y_i\partial y_o\partial y_r}G_{k,j}^v(x,y)\\
&&\phantom{\EE\left[+2\sum_{s,u,v=1}^d\right.}\cdot\left(\frac{\partial^2}{\partial y_i\partial y_o}G_{k,j}^s(x,y)\frac{\partial}{\partial y_r}G_{k,j}^u(x,y)+\frac{\partial^2}{\partial y_i\partial y_r}G_{k,j}^s(x,y)\frac{\partial}{\partial y_o}G_{k,j}^u(x,y)\right.\\
&&\phantom{\EE\left[+2\sum_{s,u,v=1}^d\cdot\left(\right.\right.}\left.+\frac{\partial}{\partial y_i}G_{k,j}^s(x,y)\frac{\partial^2}{\partial y_o\partial y_r}G_{k,j}^u(x,y)\right)\\
&&\phantom{\EE\left[\right.}+2\sum_{s,u,v,w=1}^d\frac{\partial^3}{\partial x_s\partial x_u\partial x_v}\Phi_{\Delta,k+1}^{l}(G_{k,j}(x,y))\frac{\partial}{\partial x_w}\Phi_{\Delta,k+1}^{l}(G_{k,j}(x,y))\frac{\partial^3}{\partial y_i\partial y_o\partial y_r}G_{k,j}^w(x,y)\\
&&\phantom{\EE\left[+2\sum_{s,u,v,w=1}^d\right.}\cdot\frac{\partial}{\partial y_i}G_{k,j}^s(x,y)\frac{\partial}{\partial y_o}G_{k,j}^u(x,y)\frac{\partial}{\partial y_r}G_{k,j}^v(x,y)\\
&&\phantom{\EE\left[\right.}+6d^2\sum_{s,u=1}^d\left(\frac{\partial^2}{\partial x_s\partial x_u}\Phi_{\Delta,k+1}^{l}(G_{k,j}(x,y))\right)^2\left(\left(\frac{\partial^2}{\partial y_i\partial y_o}G_{k,j}^s(x,y)\frac{\partial}{\partial y_r}G_{k,j}^u(x,y)\right)^2\right.\\
&&\phantom{\EE\left[+6d^2\sum_{s,u=1}^d\left(\frac{\partial^2}{\partial x_s\partial x_u}\Phi_{\Delta,k+1}^{l}(G_{k,j}(x,y))\right)^2\left(\right.\right.}+\left(\frac{\partial^2}{\partial y_i\partial y_r}G_{k,j}^s(x,y)\frac{\partial}{\partial y_o}G_{k,j}^u(x,y)\right)^2\\
&&\phantom{\EE\left[+6d^2\sum_{s,u=1}^d\left(\frac{\partial^2}{\partial x_s\partial x_u}\Phi_{\Delta,k+1}^{l}(G_{k,j}(x,y))\right)^2\left(\right.\right.}\left.+\left(\frac{\partial}{\partial y_i}G_{k,j}^s(x,y)\frac{\partial^2}{\partial y_o\partial y_r}G_{k,j}^u(x,y)\right)^2\right)\\
&&\phantom{\EE\left[\right.}\left.+2d^3\sum_{s,u,v=1}^d\left(\frac{\partial^3}{\partial x_s\partial x_u\partial x_v}\Phi_{\Delta,k+1}^{l}(G_{k,j}(x,y))\frac{\partial}{\partial y_i}G_{k,j}^s(x,y)\frac{\partial}{\partial y_o}G_{k,j}^u(x,y)\frac{\partial}{\partial y_r}G_{k,j}^v(x,y)\right)^2\right]
\end{eqnarray*}
Denote 
\[
\zeta_{k+1}^{i,o,r,s}=\EE\left[\left(\frac{\partial^3}{\partial y_i\partial y_o\partial y_r}G_{k+1,j}^s(x,y)\right)^{2}\right]
\]
and $e_{s,u,v}\doteq e_s+e_u+e_v$, then we get, due to $3a^2b^2c^2\le a^6+b^6+c^6$ and
\begin{align*}
2\EE\left[XYZU\right]\le & 2\sqrt{\EE\left[X^2\right]}\sqrt[6]{\EE\left[Y^6\right]}\sqrt[6]{\EE\left[Z^6\right]}\sqrt[6]{\EE\left[U^6\right]}\le \EE\left[X^2\right]+\sqrt[3]{\EE\left[Y^6\right]}\sqrt[3]{\EE\left[Z^6\right]}\sqrt[3]{\EE\left[U^6\right]}\\
&\le \EE\left[X^2\right]+\frac{1}{3}\left(\EE\left[Y^6\right]+\EE\left[Z^6\right]+\EE\left[U^6\right]\right),
\end{align*}
for $k=j+1,\ldots,J-1$,
\begin{eqnarray*}
\zeta_{k+1}^{i,o,r,l} & \leq & (1+A_{2,l}\Delta)\zeta_k^{i,o,r,l}
+\sum_{s\neq l}\left\{C_{1,1,l,e_l,e_s}\Delta(\zeta_k^{i,o,r,l}+\zeta_k^{i,o,r,s})+(d-1)B_{2,l,e_s}\Delta\zeta_k^{i,o,r,s}\right\}\\
&&+\sum_{s,u,v=1}^d C_{1,1,l,e_v,e_{s,u}}\Delta\left(\zeta_{k,2}^{i,o,r,v}+\frac{1}{2}\left(\rho_{k,4}^{i,s}+\rho_{k,4}^{o,u}+\rho_{k,4}^{r,u}+\psi_{k,4}^{i,o,s}+\psi_{k,4}^{i,r,s}+\psi_{k,4}^{o,r,u}\right)\right)\\
&&+\sum_{s,u,v,w=1}^d C_{1,1,l,e_w,e_{s,u,v}}\Delta\left(\zeta_{k,2}^{i,o,r,w}+\frac{1}{3}\left(\rho_{k,6}^{i,s}+\rho_{k,6}^{o,u}+\rho_{k,6}^{r,v}\right)\right)\\
&&+3d^2\sum_{s,u=1}^dB_{2,l,e_{s,u}}\Delta\left(\rho_{k,4}^{i,s}+\rho_{k,4}^{o,u}+\rho_{4,k}^{r,u}+\psi_{k,4}^{i,o,s}+\psi_{k,4}^{i,r,s}+\psi_{k,4}^{o,r,u}\right)\\
&&+d^3\sum_{s,u,v=1}^dB_{2,l,e_{s,u,v}}\Delta\frac{1}{3}\left(\rho_{k,6}^{i,s}+\rho_{k,6}^{o,u}+\rho_{k,6}^{r,w}\right).
\end{eqnarray*}
Further, denote 
\[
\zeta_{k+1}^{i,o,r}=\sum_{l=1}^d\zeta_{k+1}^{i,o,r,l},
\]
then we get for $k=j+1,\ldots,J-1$,
\begin{eqnarray*}
\zeta_{k+1,2}^{i,o,r} & \leq & (1+A_{2}\Delta)\zeta_{k,2}^{i,o,r}
+2(d-1)C_{1,1}\Delta\zeta_{k,2}^{i,o,r}+(d-1)^2B_{2}\Delta\zeta_{k,2}^{i,o,r}\\
&&+d^3\tilde{C}_{1,1}\Delta\left(\zeta_{k,2}^{i,o,r}+\frac{1}{2}\left(\rho_{k,4}^{i}+\rho_{k,4}^{o}+\rho_{k,4}^r+\psi_{k,4}^{i,o}+\psi_{k,4}^{i,r}+\psi_{k,4}^{o,r}\right)\right)\\
&&+d^4\tilde{\tilde{C}}_{1,1}\Delta\left(\zeta_{k,2}^{i,o,r}+\frac{1}{3}\left(\rho_{k,6}^{i}+\rho_{k,6}^{o}+\rho_{k,6}^{r}\right)\right)\\
&&+3d^4\tilde{B}_2\Delta\left(\rho_{k,4}^{i}+\rho_{k,4}^{o}+\rho_{k,4}^r+\psi_{k,4}^{i,o}+\psi_{k,4}^{i,r}+\psi_{k,4}^{o,r}\right)+d^6\tilde{\tilde{B}}_2\Delta\frac{1}{3}\left(\rho_{k,6}^{i}+\rho_{k,6}^{o}+\rho_{k,6}^{r}\right).
\end{eqnarray*}
where $\tilde{\tilde{C}}_{1,1}\doteq\max\limits_{l,s,u,v,w=1,\ldots,d}C_{1,1,l,e_w,e_{s,u,v}}$ and $\tilde{\tilde{B}}_{2}\doteq\max\limits_{l,s,u,v=1,\ldots,d}B_{2,l,e_{s,u,v}}$. 
This gives us
\[
\zeta_{k+1,2}^{i,o,r}\leq (1+\kappa_{11}\Delta)\zeta_{k,2}^{i,o,r}+\kappa_{12},\quad k=j+1,\ldots,J-1
\]
for some constants $\kappa_{11},\kappa_{12}>0$,
leading to
\[
\zeta_{k,2}^{i,o,r}\leq(1+\kappa_{11}\Delta)^{k-j-1}\zeta_{j+1,2}^{i,o,r}+\kappa_{13}=\kappa_{13}
,\,k=j+1,\ldots,J-1,
\]
where $\kappa_{13}>0$ and
\[
\zeta_{j+1,2}^{i,o,r}=\sum_{s=1}^d\EE\left[\left(\frac{\partial^3}{\partial y_i\partial y_o\partial y_r}\Phi_{\Delta}^s(x,y+\xi_{j+1})\right)^{2}\right]=0.
\]
Thus, we obtain the boundednesss of
\begin{eqnarray*}
&&\frac{\partial^3}{\partial y_i\partial y_o\partial y_r}u_{\Delta}(t_{j},x,y)\\
&=&\EE\left[\sum_{s=1}^d\frac{\partial }{\partial x_s}f(G_{J,j}(x,y))\frac{\partial^3}{\partial y_i\partial y_o\partial y_r}G_{J,j}^s(x,y)\right.\\
&&\phantom{\EE\left[\right.}+\sum_{s,u=1}^d\frac{\partial^2 }{\partial x_s\partial x_u}f(G_{J,j}(x,y))\left(\frac{\partial^2}{\partial y_i\partial y_o}G_{J,j}^s(x,y)\frac{\partial}{\partial y_r}G_{J,j}^u(x,y)+\frac{\partial^2}{\partial y_i\partial y_r}G_{J,j}^s(x,y)\frac{\partial}{\partial y_o}G_{J,j}^u(x,y)\right.\\
&&\phantom{\EE\left[+\sum_{s,u=1}^d\frac{\partial^2 }{\partial x_s\partial x_u}f(G_{J,j}(x,y))\left(\right.\right.}\left.+\frac{\partial}{\partial y_i}G_{J,j}^s(x,y)\frac{\partial^2}{\partial y_o\partial y_r}G_{J,j}^u(x,y)\right)\\
&&\phantom{\EE\left[\right.}\left.+\sum_{s,u,v=1}^d\frac{\partial^3 }{\partial x_s\partial x_u\partial x_v}f(G_{J,j}(x,y))\frac{\partial}{\partial y_i}G_{J,j}^s(x,y)\frac{\partial}{\partial y_o}G_{J,j}^u(x,y)\frac{\partial}{\partial y_r}G_{J,j}^v(x,y)\right],
\end{eqnarray*}
provided that $\sigma_{ki}$ and all the derivatives of order 1, 2 and 3 for $f,\mu_k,\sigma_{ki}$ are bounded.

\subsection*{Proof of Lemma~\protect\ref{lemm_cv_strong}}
Cf.~Theorem~5.2 in~\cite{belomestny2016variance}.

\subsection*{Proof of Theorem~\protect\ref{theo_cv_strong}}
We have,
by the martingale property of
$(\tilde M^{{int},1}_{\Delta,j\Delta})_{j=0,\ldots,J}$,
where $\tilde M^{{int},1}_{\Delta,j\Delta}$
is given by~\eqref{eq:cv_strong2} with $J$
being replaced by~$j$,
and by the orthogonality of the system
$\Delta_j W^i$,
\begin{eqnarray*}
\Var[f(X_{\Delta,T})-
\tilde{M}^{{int},1}_{\Delta,T}]&=&\Var[f(X_{\Delta,T})-
M^{{int},1}_{\Delta,T}]+\Var[M^{{int},1}_{\Delta,T}-
\tilde{M}^{{int},1}_{\Delta,T}]\\
&\lesssim&\frac{1}{J}+\Delta\sum_{j=1}^J
\sum_{i=1}^m\EE\|\sum_{k=1}^d
(\tilde g_{j,k}-g_{j,k})\sigma_{ki}\|^2_{L^2(\PP_{\Delta,j-1})}
\\
&\le&\frac{1}{J}+
d\Delta\sum_{j=1}^J
\sum_{i=1}^m\sum_{k=1}^d
\EE\|(\tilde g_{j,k}-g_{j,k})\sigma_{ki}\|^2_{L^2(\PP_{\Delta,j-1})}\\
&\le&
\frac{1}{J}+d^2Tm
\sigma_{\max} ^2
\left\{\tilde c\left(\Sigma+A^2(\log N+1)\right)\frac{\binom{p+d}d Q^d}{N}\right.
\\
&&
\left.+\frac{8\,C_h^2}{(p+1)!^2 d^{2-2/h}}
\left(\frac{Rd}Q\right)^{2p+2}
+8A^2 B_\nu R^{-\nu}\right\}.
\end{eqnarray*}

\subsection*{Proof of Theorem~\protect\ref{theo_comp_strong}}
Let us, for simplicity, first ignore the $\log(N)$-term in~\eqref{mse_strong_constr} and only consider the terms w.r.t. the variables $J,N,N_0,Q,R$ which shall be optimised, since the constants $d,m,c_{p,d}$, $(p+1)!$, $B_\nu$ do not affect the terms on $\varepsilon$. Further, we consider the log-cost and log-constraints rather than~\eqref{cost_strong} and ~\eqref{mse_strong_constr}. Let us subdivide the optimisation problem into two cases:
\begin{enumerate}
\item $N\lesssim N_0$.
This gives us the Lagrange function
\begin{align}
\label{lagrange_strong}
&L_{\lambda_1,\ldots,\lambda_6}(J,N,N_0,Q,R)\\
\notag
\doteq & \log(J)+\log(N_0)+d\log(Q)+\lambda_1(-2\log(J)-2\log(\varepsilon))\\
\notag
&+\lambda_2(-\log(J)-\log(N_0)-2\log(\varepsilon))\\
\notag
&+\lambda_3(d\log(Q)-\log(N)
-\log(N_0)-2\log(\varepsilon))\\
\notag
&+\lambda_4(2(p+1)(\log(R)-\log(Q))-\log(N_0)-2\log(\varepsilon))\\
\notag
&+\lambda_5(-\nu\log(R)-\log(N_0)
-2\log(\varepsilon))+\lambda_6(\log(N)-\log(N_0)),
\end{align}
where $\lambda_1,\ldots,\lambda_6\ge 0$. Thus, considering of the conditions $\frac{\partial L}{\partial J}=\frac{\partial L}{\partial N}=\frac{\partial L}{\partial N_0}=\frac{\partial L}{\partial Q}=\frac{\partial L}{\partial R}\stackrel{!}{=}0$ gives us the following relations
\begin{align*}
\lambda_1&=\frac{1-\lambda_2}{2},\\
\lambda_3&=\frac{2(p+1)(\nu(1-\lambda_2)-d)-d\nu}{d\nu+2(p+1)(d+2\nu)}=\lambda_6,\\
\lambda_4&=\frac{d\nu(3-\lambda_2)}{d\nu+2(p+1)(d+2\nu)},\\
\lambda_5&=\frac{2d(p+1)(3-\lambda_2)}{d\nu+2(p+1)(d+2\nu)}.
\end{align*}
The case $\lambda_1,\ldots,\lambda_6>0$ is not feasible, since all constraints in~\eqref{lagrange_strong} can not be active, that is they cannot become zero simultaneously because of six (linearly independent) equalities on five unknowns. Hence, we derive the solutions under $\lambda_i=0$ for different $i$ and observe which one is actually optimal.
\begin{enumerate}
\item $\lambda_1=0\;\Rightarrow\;\lambda_3=\lambda_6=-\frac{d(2(p+1)+\nu)}{d\nu+2(p+1)(d+2\nu)}<0$. Due to negative $\lambda_3,\lambda_6$, this case is not optimal.
\item $\lambda_2=0\,\Rightarrow\;\lambda_1,\lambda_4,\lambda_5>0,\lambda_3=\lambda_6=\frac{2(p+1)(\nu-d)-d\nu}{d\nu+2(p+1)(d+2\nu)}$. Again, we make a case distinction:
\begin{enumerate}
\item $\lambda_3=\lambda_6=0\;\Rightarrow\;\nu=\frac{2d(p+1)}{2(p+1)-d}$ for $2(p+1)>d$. This gives us, due to $\lambda_1,\lambda_4,\lambda_5>0$,
\begin{align*}
J&\asymp\varepsilon^{-1},\\
Q&\asymp\left[\frac{1}{N_0\varepsilon^2}\right]^\frac{1}{d},\\
JQ^dN_0&\asymp \varepsilon^{-3}.
\end{align*}
This solution is no improvement compared to the SMC approach.
\item $\lambda_3=\lambda_6>0\;\Rightarrow\;\nu>\frac{2d(p+1)}{2(p+1)-d}$ for $2(p+1)>d$. In this case, all constraints apart from the second one in~\eqref{lagrange_strong}, corresponding to $\lambda_2$, are active. Then we obtain
\begin{align*}
J&\asymp \varepsilon^{-1},\\
Q&\asymp \varepsilon^{-\frac{2\nu+4(p+1)}{d\nu+2(p+1)(d+2\nu)}}, \\
N_0&\asymp \varepsilon^{-\frac{2d\nu+4(p+1)(d+\nu)}{d\nu+2(p+1)(d+2\nu)}}, \\
JQ^dN_0&\asymp \varepsilon^{-\frac{5d\nu+2(p+1)(5d+4\nu)}{d\nu+2(p+1)(d+2\nu)}},
\end{align*}
which is a better solution than the previous one. Moreover, the remaining constraint $\frac{1}{JN_0}\lesssim\varepsilon^2$ is also satisfied under this solution.
\end{enumerate}
\item $\lambda_3=\lambda_6=0\,\Rightarrow\;\lambda_1,\lambda_4,\lambda_5>0,\lambda_2=\frac{2(p+1)(\nu-d)-d\nu}{2(p+1)\nu}$. The case $\lambda_2=0$ is the same as the last but one and thus gives us $JQ^d N_0\asymp \varepsilon^{-3}$. The case $\lambda_2>0$ leads to four active constraints in~\eqref{lagrange_strong}, namely the ones corresponding to $\lambda_1,\lambda_2,\lambda_4,\lambda_5$, such that
\begin{align*}
J&\asymp \varepsilon^{-1},\\
Q&\asymp \varepsilon^{-\frac{\nu+2(p+1)}{2\nu(p+1)}}, \\
N_0&\asymp \varepsilon^{-1}, \\
JQ^dN_0&\asymp \varepsilon^{-\frac{d\nu+2(p+1)(d+2\nu)}{2\nu(p+1)}}.
\end{align*}
This solution seems to be nice at the first moment. However, it does not satisfy both  constraints corresponding to $\lambda_3,\lambda_6$. On the one hand, we have for the third constraint $N\gtrsim\varepsilon^{-1-\frac{d\nu+2d(p+1)}{2\nu(p+1)}}$. On the other hand, we have for the sixth constraint $N\lesssim\varepsilon^{-1}$. Hence, this is not an admissible solution.
\item $\lambda_4=0\,\Rightarrow\;\lambda_1=-1$. Since $\lambda_1$ is negative, this case is not optimal.
\item $\lambda_5=0\,\Rightarrow\;\lambda_1=-1$. As for the previous one, this case is not optimal.
\end{enumerate}
\item $N\gtrsim N_0$.
This gives us the Lagrange function
\begin{align*}
&
\tilde{L}_{\lambda_1,\ldots,\lambda_6}(J,N,N_0,Q,R)\\
\doteq & \log(J)+\log(N)+d\log(Q)+\lambda_1(-2\log(J)-2\log(\varepsilon))\\
&+\lambda_2(-\log(J)-\log(N_0)-2\log(\varepsilon))\\
&+\lambda_3(d\log(Q)-\log(N)
-\log(N_0)-2\log(\varepsilon))\\
&+\lambda_4(2(p+1)(\log(R)-\log(Q))-\log(N_0)-2\log(\varepsilon))\\
&+\lambda_5(-\nu\log(R)-\log(N_0)
-2\log(\varepsilon))+\lambda_6(\log(N_0)-\log(N)).
\end{align*}
Analogously to the procedure above we get the same optimal solution, that is
\begin{align*}
J&\asymp \varepsilon^{-1},\\
Q&\asymp \varepsilon^{-\frac{2\nu+4(p+1)}{d\nu+2(p+1)(d+2\nu)}}, \\
N&\asymp \varepsilon^{-\frac{2d\nu+4(p+1)(d+\nu)}{d\nu+2(p+1)(d+2\nu)}}, \\
JQ^dN&\asymp \varepsilon^{-\frac{5d\nu+2(p+1)(5d+4\nu)}{d\nu+2(p+1)(d+2\nu)}}.
\end{align*}
\end{enumerate}
Now we consider also the remaining terms $c_{p,d}$, $(p+1)!$, $B_\nu$ and obtain~\eqref{J_strong}--\eqref{compl_strong} via equalising all constraints in~\eqref{mse_strong_constr} apart from the second one. Finally, we add the $\log$-term concerning $\varepsilon$ in the parameters $N,N_0$ to ensure that all constraints are really satisfied.

\subsection*{Proof of Lemma~\protect\ref{lemm_cv_strong_ser}}
Cf.~Theorem 5.2 in~\cite{belomestny2016variance}.

\subsection*{Proof of Theorem~\protect\ref{theo_cv_strong_ser}}
We have,
by the martingale property of
$(\tilde M^{{ser},1}_{\Delta,j\Delta})_{j=0,\ldots,J}$,
where $\tilde M^{{ser},1}_{\Delta,j\Delta}$
is given by~\eqref{eq:cv_strong_ser} with $J$
being replaced by~$j$,
and by the orthonormality of the system
$\frac{\Delta_j W}{\sqrt{\Delta}}$,
\begin{eqnarray*}
\Var[f(X_{\Delta,T})-
\tilde{M}^{{ser},1}_{\Delta,T}]&=&\Var[f(X_{\Delta,T})-
M^{{ser},1}_{\Delta,T}]+\Var[M^{{ser},1}_{\Delta,T}-
\tilde{M}^{{ser},1}_{\Delta,T}]\\
&\lesssim&\frac{1}{J}+\sum_{j=1}^J\sum_{i=1}^m
\EE\|
\tilde a_{j,e_i}-a_{j,e_i}\|^2_{L^2(\PP_{\Delta,j-1})}
\\
&\le&
\frac{1}{J}+Jm
\left\{\tilde c\left(\Sigma+A^2\Delta(\log N+1)\right)\frac{c_{p,d} Q^d}{N}\right.
\\
&&
\left.+\frac{8\,C_h^2}{(p+1)!^2}
\left(\frac{R}Q\right)^{2p+2}
+8A^2\Delta B_\nu R^{-\nu}\right\}.
\end{eqnarray*}

\subsection*{Proof of Theorem~\protect\ref{theo_comp_strong_ser}}
The proof is similar to the one of Theorem~\ref{theo_comp_strong}.

\bibliographystyle{abbrv}
\bibliography{bibliography}

\begin{thebibliography}{10}

\bibitem{AAO}
J.~Akahori, T.~Amaba, and K.~Okuma.
\newblock A discrete-time {C}lark-{O}cone formula and its application to an
  error analysis.
\newblock {\em Preprint, arXiv:1307.0673v2}, 2013.

\bibitem{belomestny2016variance}
D.~Belomestny, S.~H\"afner, T.~Nagapetyan, and M.~Urusov.
\newblock Variance reduction for discretised diffusions via regression.
\newblock {\em Preprint, arXiv:1510.03141v3}, 2016.

\bibitem{giles2008multilevel}
M.~B. Giles.
\newblock Multilevel {M}onte {C}arlo path simulation.
\newblock {\em Operations Research}, 56(3):607--617, 2008.

\bibitem{gyorfi2002distribution}
L.~Gy{\"o}rfi, M.~Kohler, A.~Krzy{\.z}ak, and H.~Walk.
\newblock {\em A distribution-free theory of nonparametric regression}.
\newblock Springer Series in Statistics. Springer-Verlag, New York, 2002.

\bibitem{karatzas2012brownian}
I.~Karatzas and S.~Shreve.
\newblock {\em Brownian motion and stochastic calculus}, volume 113.
\newblock Springer Science \& Business Media, 2012.

\bibitem{KP}
P.~Kloeden and E.~Platen.
\newblock {\em Numerical solution of stochastic differential equations},
  volume~23.
\newblock Springer, 1992.

\bibitem{MilsteinTretyakov:2004}
G.~N. Milstein and M.~V. Tretyakov.
\newblock {\em Stochastic numerics for mathematical physics}.
\newblock Scientific Computation. Springer-Verlag, Berlin, 2004.

\bibitem{milstein2009practical}
G.~N. Milstein and M.~V. Tretyakov.
\newblock Practical variance reduction via regression for simulating
  diffusions.
\newblock {\em SIAM Journal on Numerical Analysis}, 47(2):887--910, 2009.

\bibitem{muller2015complexity}
T.~M{\"u}ller-Gronbach, K.~Ritter, and L.~Yaroslavtseva.
\newblock On the complexity of computing quadrature formulas for marginal
  distributions of sdes.
\newblock {\em Journal of Complexity}, 31(1):110--145, 2015.

\bibitem{newton1994variance}
N.~J. Newton.
\newblock Variance reduction for simulated diffusions.
\newblock {\em SIAM journal on applied mathematics}, 54(6):1780--1805, 1994.

\end{thebibliography}

\end{document}